%% file: main_bc.tex
\documentclass[a4paper,10pt]{article}
\usepackage[top=13mm, bottom=23mm, inner=20mm, outer=20mm]{geometry}

\input{Templates/StandardPackagesAndFormatting}
\input{Templates/StandardMathsMacros}
\input{Templates/StandardTheorems}
\usepackage{natbib,doi}
\graphicspath{{Figures/}}

\usepackage[scr=boondox]{mathalpha}

\DeclareMathOperator{\tri}{tri}
\DeclareMathOperator{\sqr}{sqr}

\newcommand{\whatchar}[1]{\widehat{\,#1\,}}
\DeclareMathOperator\TV{TV}

\begin{document}

\input{Sections/FrontMatter}

\input{Sections/Introduction}
\input{Sections/Discretization}
\input{Sections/Spatial}

\input{Sections/Temporal}
\input{Sections/Test}
\input{Sections/Discussion}
\input{Sections/Acknowledgements}

\begin{myappendices}
	\input{Sections/RK_contraction}
	\input{Sections/Implementation}
\end{myappendices}

\bibliographystyle{plainnat}
\bibliography{Bibliography/Books,Bibliography/Mine,Bibliography/NumericalMethodsHyperbolic,Bibliography/ShallowWater,Bibliography/AnalysisHyperbolic,Bibliography/GravityCurrents}

\end{document}

%% file: Templates/StandardPackagesAndFormatting.tex
\usepackage{mathtools,amssymb,mathrsfs,bm}

\usepackage{graphicx}
\usepackage{caption,subcaption,float}
\usepackage{tabularx,multirow}
\usepackage[inline]{enumitem}

\usepackage{xcolor}

\numberwithin{equation}{section}
\numberwithin{figure}{section}
\numberwithin{table}{section}

\usepackage{calc,etoolbox}

\usepackage[inline]{enumitem}

\makeatletter
\@ifpackageloaded{hyperref}{}{\usepackage[hidelinks]{hyperref}}
\@ifpackageloaded{cleveref}{}{
	\usepackage[compress]{cleveref}

	\crefformat{section}{\S#2#1#3}
	\crefformat{subsection}{\S#2#1#3}
	\crefformat{subsubsection}{\S#2#1#3}
	\crefrangeformat{section}{\S\S#3#1#4 to~#5#2#6}
	\crefmultiformat{section}{\S\S#2#1#3}{ and~#2#1#3}{, #2#1#3}{ and~#2#1#3}
	
	\crefname{figure}{fig.}{figs.}
	\Crefname{figure}{Fig.}{Figs.}
	
	\crefname{equation}{}{}
	\crefformat{equation}{(#2#1#3)}
	\crefmultiformat{equation}{(#2#1#3)}{ and (#2#1#3)}{, (#2#1#3)}{ and (#2#1#3)}
}
\makeatother

\usepackage[title]{appendix}
\newenvironment{myappendices}
{
	\begin{appendices}
	\crefalias{section}{appsec}
}{
	\end{appendices}
}

\crefname{appsec}{Appendix}{Appendices}
\Crefname{appsec}{Appendix}{Appendices}

\usepackage{xspace}
\makeatletter
\DeclareRobustCommand\onedot{\futurelet\@let@token\@onedot}
\def\@onedot{\ifx\@let@token.\else.\null\fi\xspace}
\makeatother
\def\eg{e.g\onedot} 
\def\ie{i.e\onedot} 
\def\cf{\textit{cf}\onedot}

\def\sf{s.f\onedot}

\makeatletter
\newcommand{\floattag}[1]{%
  \@namedef{the\@captype}{#1}%
  \@namedef{theH\@captype}{#1}%
  \addtocounter{\@captype}{-1}}
\makeatother

%% file: Templates/StandardTheorems.tex
\usepackage{amsthm, thmtools, thm-restate}
\theoremstyle{plain}
\declaretheorem[numberwithin=section,	refname={theorem,theorems},			Refname={Theorem,Thorems}]								{theorem}
\declaretheorem[name=Lemma,				refname={lemma,lemmas},				Refname={Lemma,Lemmas},				sibling=theorem]	{lemma}
\declaretheorem[name=Proposition,		refname={proposition,propositions},	Refname={Proposition,Propositions},	sibling=theorem]	{prop}

\theoremstyle{definition}

\declaretheorem[name=Discussion,		refname={discussion,discussions},	Refname={Discussion,Discussions},	numbered=no]		{discuss}
\theoremstyle{remark}

\theoremstyle{definition}
\declaretheorem[name=Test problem, numberwithin=section]{testproblem}

\crefformat{testproblem}{test problem #2#1#3}
\crefrangeformat{testproblem}{test problems #3#1#4 to~#5#2#6}
\Crefrangeformat{testproblem}{Test problems #3#1#4 to~#5#2#6}
\crefmultiformat{testproblem}{test problems #2#1#3}{ and~#2#1#3}{, #2#1#3}{ and~#2#1#3}

%% file: Sections/FrontMatter.tex
\title{The implementation of a broad class of boundary conditions for non-linear hyperbolic systems}
\author{Edward W. G. Skevington}
\date{21/06/2021}
\maketitle

\begin{abstract}
	We propose methods that augment existing numerical schemes for the simulation of hyperbolic balance laws with Dirichlet boundary conditions to allow for the simulation of a broad class of differential algebraic conditions. Our approach is similar to that of Thompson (1987), where the boundary values were simulated by combining characteristic equations with the time derivative of the algebraic conditions, but differs in two important regards. Firstly, when the boundary is a characteristic of one of the fields Thompson's method can fail to produce reasonable values. We propose a method of combining the characteristic equations with extrapolation which ensures convergence. Secondly, the application of algebraic conditions can suffer from $\order{1}$ drift-off error, and we discuss projective time-stepping algorithms designed to converge for this type of system. Test problems for the shallow water equations are presented to demonstrate the result of simulating with and without the modifications discussed, illustrating their necessity for certain problems.
\end{abstract}

%% file: Sections/Introduction.tex
\section{Introduction}	\label{sec:Intro}

An important class of equations in the study of fluid dynamics, in areas such as gas dynamics, hydraulics, and gravity currents, is non-linear hyperbolic balance laws,
\begin{align}	\label{eqn:hyp_system}
	\pdv{Q}{t} + \pdv{}{x} \ppar{F} &= \Psi
\end{align}
where $Q:\ppar{x,t} \mapsto \mathbb{R}^M$ is a vector of densities, which are transported by the flux $F:\ppar{Q,x,t} \mapsto \mathbb{R}^M$ and created/destroyed by an algebraic source term $\Psi:\ppar{Q,x,t} \mapsto \mathbb{R}^M$. (See \eg \cite{bk_Serre_SCL1,bk_Serre_SCL2} for analytic properties, and \cite{bk_Leveque_FVM} for numerical properties.) Here we only consider one spatial dimension, but the ideas presented are generalisable to multidimensional problems. For \cref{eqn:hyp_system} to be classed as hyperbolic, we require that $A \eqdef \pdv*{F}{Q}$ is diagonalisable with real eigenvalues $\lambda^{(m)}$ ordered as $\lambda^{(1)} \leq \ldots \leq \lambda^{(M)}$ with corresponding left and right eigenvectors $l^{(m)}$ and $r^{(m)}$ respectively. Letting $L$ and $R$ be $M \times M$ matrices ($L,R \in \mathbb{M}\ppar{M,M}$) with elements $L_{ij} \eqdef l_{j}^{(i)}$, $R_{ij} \eqdef r_{i}^{(j)}$ where $LR = I$ we have $LAR = \Lambda \eqdef \diag \pbrk*{\lambda^{(1)} \ldots \lambda^{(M)}}$. In regions where the solution is continuous \cref{eqn:hyp_system} can be written in characteristic form
\begin{align}	\label{eqn:hyp_characteristics}
	L\pdv{Q}{t} + \Lambda L \pdv{Q}{x} &= L \Psi'
	\qquad\text{where}\qquad
	\Psi' \eqdef \Psi - \pdv{F}{x}.
\end{align}

A problem on a finite domain $x_L\ppar{t} \leq x \leq x_R\ppar{t}$ requires boundary conditions. Throughout this article we will discuss the right boundary $x=x_R$ at which $Q=Q_R$, the left treated similarly by symmetry. For the case of a stationary boundary point ($x_R$ constant) it has been shown that the imposition of an inhomogeneous Dirichlet boundary condition is sufficient to close the system \cite{ar_Nordstrom_2005,ar_Guaily_2012}, that is
\begin{align}	\label{eqn:BC_dirichlet}
	L^- Q_R &= L^- Q_{R0}
	\qquad\text{where}\qquad
	L^- \eqdef
	\pbrk*{\begin{array}{@{}c@{}c@{}c@{}}
		\leftarrow	& l^{(1)}	& \rightarrow		\\
				 	& \vdots	& 					\\
		\leftarrow	& l^{(M_I)}	& \rightarrow	
	\end{array}},
\end{align}
$M_I$ is the largest value for which $\lambda^{(M_I)} < 0$, and $Q_{R0}$ some function of $t$. If $x_R$ varies in time then, by a straightforward change of variables, \cref{eqn:BC_dirichlet} with $M_I$ the largest value so that $\lambda^{(M_I)} < \dot{x}_R \eqdef \dv*{x_R}{t}$ is sufficient to close the system. This imposes one condition on each of the $M_I$ \emph{incoming} characteristics for which $\lambda^{(m)} < \dot{x}_R$, and no condition on the \emph{static} and \emph{outgoing} characteristics for which $\lambda^{(m)} = \dot{x}_R$ and $\lambda^{(m)} > \dot{x}_R$ respectively. All that remains is to specify $L^- Q_{R0}$ and the location of the boundary $x_R$ for the problem at hand, for a total of $M_B \eqdef M_I+1$ conditions.

Condition \cref{eqn:BC_dirichlet} covers many important cases where the value of $L^- Q_{R0}$ is readily deduced from the boundary conditions, for example the slip condition for the Euler equations \cite{ar_Svard_2014b} (and the no-slip condition for the parabolic Navier-Stokes equations, \cite{ar_Svard_2008}). We stress that, while a condition of the type \cref{eqn:BC_dirichlet} is \emph{sufficient} to close the system, it is by no means \emph{necessary}. For example, in \cite{ar_Skevington_F001_Draining} the shallow water equations were examined for the partial collapse of a dam. Denoting the depth of the fluid $h$ and the volume flux $q$ the dimensionless boundary condition used was
\begin{align}\label{eqn:BC_example}
	\frac{q^2}{h^2} + 2 \ppar*{h-h_b} - 3 q^{2/3} &= 0
	&&\text{at}&
	x &= x_R
\end{align}
with $0 \leq h_b \leq h$, and $0 \leq q \leq h^{3/2}$, where $h_b$ is the elevation of the barrier resulting from the collapse. This boundary condition is not of the form \cref{eqn:BC_dirichlet}, and yet yielded a unique solution. Many important boundary conditions for physical problems (including \cref{eqn:BC_dirichlet,eqn:BC_example}) belong to the differential-algebraic (DA) class, for which
\begin{subequations}\label{eqn:BC_general}	\begin{align}
	B_Q \ppar*{Q_R,Q_{xR},\dot{x}_R,x_R,t} \cdot \dv{Q_R}{t} + B_x \ppar*{Q_R,Q_{xR},\dot{x}_R,x_R,t} \cdot \dv[2]{x_R}{t} &= b \ppar*{Q_R,Q_{xR},\dot{x}_R,x_R,t},
	\label{eqn:BC_general_diff}	\\
	g \ppar*{Q_R,Q_{xR},\dot{x}_R,x_R,t} &=0
	\label{eqn:BC_general_alge}
	\qquad\text{where}\qquad
	Q_{xR} = \eval*{\pdv{Q}{x}}_{\mathrlap{x=x_R}}.
\end{align}\end{subequations}
Here $\ppar*{Q_R^T , \dot{x}_R}^T$ is the data to be established at the boundary, where we include the velocity of the boundary as an unknown to be solved for. Each of $B_Q$, $B_x$, and $b$ has $M_D$ rows to enforce $M_D$ differential boundary conditions, while $g$ has $M_A$ rows to enforce $M_A$ algebraic conditions. We require $M_D + M_A = M_B$, to enforce the correct number of conditions. The general system \cref{eqn:hyp_system,eqn:BC_general} is so broad that it will contain some subset of ill-posed problems; in what follows we will assume that the system of interest at least satisfies existence and uniqueness.

A final special type of condition is the non-reflecting condition \cite{ar_Hedstrom_1978}. This is employed in numerical schemes when only a portion of the physical domain is simulated, and attempts to enforce that no waves enter the domain. Thus for numerical schemes we may wish to implement non-reflecting condition, a Dirichlet condition, or a DA condition. We focus on the latter, but also discuss non-reflecting conditions.

Boundary conditions for numerical schemes are enforced in a variety of methods. Finite volume schemes often employ the ghost cell technique (\eg \cite{bk_Leveque_FVM}), were the boundary condition is enforced by imposing values in cells exterior to the physical domain. Specialist techniques exist for the imposition of the non-reflecting condition \cite{ar_Gross_2007}. For Dirichlet conditions it may be possible to combine an extrapolation from the bulk values with \cref{eqn:BC_dirichlet} to obtain $Q_{R0}$, and thereby establish ghost cell values \cite{ar_Du_2018}. A similar situation exists for summation-by-parts schemes introduced in \cite{ar_Strand_1994}, see \cite{ar_Fernandez_2014,ar_Svard_2014} for reviews, and \cite{ar_Fisher_2011,ar_Gassner_2016} for example schemes. Direct implementations of non-reflecting and Dirichlet conditions are available, for which convergence can be proved.

However, a direct implementation of \cref{eqn:BC_general} does not exist for finite volume, summation-by-parts, or any other scheme. Considering a semi-discrete scheme at some instance of time, the value value of $Q_R$ may be imposed on the scheme in the bulk using some standard method for Dirichlet conditions, thus the computation of $\dv*{Q_R}{t}$ and $\dv*{\dot{x}_R}{t}$ is all that remains to close the discrete system, which is what we focus on here. Following the discussion in \cite{ar_Thompson_1990}, we combine the boundary conditions \cref{eqn:BC_general} with the equations for the outgoing characteristics
\begin{align}
	L^+\pdv{Q}{t} + \Lambda^+ L^+ \pdv{Q}{x} &= L^+ \Psi'
	&&\text{where}&
	\Lambda^+ &\eqdef 
	\pbrk*{\begin{array}{@{}c@{}c@{}c@{}}
		\lambda^{(M_I+1)}	&			& 					\\
					 	& \ddots	& 					\\
						&			& \lambda^{(M)}	
	\end{array}},
	&
	L^+ &\eqdef
	\pbrk*{\begin{array}{@{}c@{}c@{}c@{}}
		\leftarrow	& l^{(M_I+1)}	& \rightarrow		\\
				 	& \vdots	& 					\\
		\leftarrow	& l^{(M)}	& \rightarrow	
	\end{array}}.
\end{align}
The approach employed in \cite{ar_Thompson_1987,ar_Thompson_1990,ar_Kim_2000,ar_Kim_2004}, which focus on the slip and no slip conditions, is to enforce the time derivative of the boundary condition $\dv*{g}{t} = 0$. Using this we obtain the system of equations
\begin{align}\label{eqn:BC_general_Thompson}
	\begin{bmatrix}
		L^+				&	0			\\
		B_Q				&	B_x			\\
		\pdv*{g}{Q_R}	&	\pdv*{g}{\dot{x}_R}
	\end{bmatrix}
	\cdot
	\dv{}{t}
	\begin{bmatrix}
		Q_R	\\
		\dot{x}_R
	\end{bmatrix}
	&=
	\begin{bmatrix}
		(\dot{x}_R I-\Lambda^+ ) L^+ Q_{xR} + L^+ \Psi'	\\
		b												\\
		- \pdv*{g}{Q_{xR}} \cdot \pdv*{Q_{xR}}{t} - \pdv*{g}{x_R} \cdot \dot{x}_R - \pdv*{g}{t}
	\end{bmatrix},
\end{align}
which are approximated by the numerical scheme.

The approach \cref{eqn:BC_general_Thompson} has a number of technical difficulties, to which we propose resolutions. The first problem occurs with static characteristic fields, for which the characteristic equations enforce
\begin{equation}	\label{eqn:hyp_characteristic_static}
	l^{(m)} \dv{Q_R}{t} = l^{(m)} \Psi'.
\end{equation}
The solution to this equation is independent of the solution in the bulk, which is qualitatively distinct from the cases of incoming and outgoing characteristics which are coupled to the solution in the bulk by the approximation to $\pdv*{Q}{x}$. Therefore, in the high resolution limit $\Delta x \rightarrow 0$, the solution can fail to satisfy
\begin{equation} \label{eqn:boundary_limit}
	Q_R(t) = \lim_{t' \rightarrow t} \lim_{x \rightarrow x_R^-} Q\ppar*{x,t'}.
\end{equation}
However, it is reasonable to insist that \cref{eqn:boundary_limit} is satisfied whenever the limit exists (the only time it will not is when a shock enters/leaves the domain), indeed \cref{eqn:boundary_limit} could be viewed as our definition of $Q_R$. The discrepancy is only important when some value of $L^- Q_R$ is dependent on the value carried by the static characteristic, which may well be the case for members of the DA class of conditions. To enforce \cref{eqn:boundary_limit} we combine the evolution equation \cref{eqn:hyp_characteristic_static} with extrapolation from the bulk which couples the dynamics at the boundary to those in the bulk. A large portion of the new material for this paper constitutes a detailed discussion of how to include extrapolation in the boundary conditions to produce convergent simulations.

The second problem is due to enforcing $\dv*{g}{t} = 0$. Any error $g \neq 0$ is carried forward in time, and numerical errors compound. This $\order{1}$ drift-off error is characteristic of approximating a system of differential-algebraic equations (DAEs) using standard time integration techniques \cite{bk_Hairer_SDAP}. The study of DAEs has produced a wide range of implicit Runge-Kutta (RK) methods that can be applied directly to systems such as ours (\eg \cite{bk_Kunkel_DAE}). However, these do not have the desirable properties for hyperbolic systems that have been proven for explicit RK methods (\eg the total variation diminishing (TVD) schemes of \cite{ar_Shu_1988}), and applying an implicit scheme across the entire domain would be computationally inefficient. We discuss the use of schemes which project onto the manifold $g=0$, and give details of the scheme from \cite{ar_Skevington_F001_Draining} which permits implicit resolution of $g=0$ while being explicit in the bulk.

The paper is organised as follows. In \cref{sec:discrete_time} we overview discretization, and present a result on the convergence of the Runge-Kutta schemes from \cite{ar_Shu_1988}. We then introduce extrapolation in \cref{sec:Extrapolation}, and assemble the system of equations to solve at the boundary for both DA and non-reflecting conditions. Spacial discretization of the boundary system is performed in \cref{sec:bc_discrete} ensuring convergence on values consistent with the bulk, and stability results are presented bounding the coefficient. The natural transformation between extrapolation and diffusion is presented in \cref{sec:Diffusion}, providing a justification of our approach in terms of the vanishing viscosity limit. \Cref{sec:TDP} is a detailed discussion of temporal discretisation and projection (\cref{sec:TDP_local,sec:TDP_projection}), including results for the scheme from \cite{ar_Skevington_F001_Draining} (\cref{sec:TDP_RKNR}), as well as a selection of numerical tests on simple problems (\cref{sec:TDP_tests}). We perform tests of the full numerical scheme using shallow water equations in \cref{sec:FVTests}, transforming onto a fixed domain in \cref{sec:TED,sec:SW}, assembling the finite volume scheme in \cref{sec:FV}, and the results of our tests are presented in \cref{sec:TP}. We conclude in \cref{sec:Discussion}.

%% file: Sections/Discretization.tex
\section{Discretization of the bulk} \label{sec:discrete_time}

\begin{table}[tp!]
	\centering
	\begin{tabular}{l c c c c}
		Order	&	$S$	&	$\eta^{[0]}$	&	$\eta^{[1]}$	&	$\eta^{[2]}$	\\
		\hline
		$1$		&	$1$	&	$0$			&	---			&	---			\\
		$2$		&	$2$	&	$0$			&	$1/2$		&	---			\\
		$3$		&	$3$	&	$0$			&	$3/4$		&	$1/3$
	\end{tabular}
	\caption{Coefficients for the convex Euler Runge-Kutta methods from \cite{ar_Shu_1988}.}
	\label{tab:RKcoeff}
\end{table}

To spatially discretize we suppose that at time $t$ the bulk solution is approximated at $J$ points $x_j(t)$ where $x_L < x_1 < x_2 < \ldots < x_J < x_R$ at which the approximated solution is $Q_j(t)$. The grid spacing may be non-uniform and we denote $\Delta x_{1/2} \eqdef x_1 - x_L$, $\Delta x_{j+1/2} \eqdef x_{j+1} - x_{j}$, and $\Delta x_{J+1/2} \eqdef x_{R} - x_{J}$. For example, a (second order) finite volume scheme has cell interfaces at $x_L = x_{1/2} < x_{3/2} < \ldots < x_{J+1/2} = x_R$ so that each cell is of width $\Delta x_j = x_{j+1/2}-x_{j-1/2}$, the approximation points are the cell centres $x_{j} = \ppar*{x_{j-1/2} + x_{j+1/2}}/2$, and the values $Q_j$ are the cell averaged values.

From the bulk scheme $\dv*{Q_{j}}{t}$ are known, and $\dv*{Q_L}{t}$, $\dv*{\dot{x}_L}{t}$, $\dv*{Q_R}{t}$, and $\dv*{\dot{x}_R}{t}$ are found using the approach we develop in \cref{sec:bc_discrete}. To evolve in time we employ the convex Euler RK schemes from \cite{ar_Shu_1988} for all differential equations; these schemes are of interest because they permit the extrapolation of properties (like TVD) from Euler time-stepping to higher order. We discretize time as $t^0<t^1<t^2<\ldots$, $\Delta t^n = t^{n+1}-t^{n}$, and denote the solution at time $t^n$ by $\pbrc*{Q}^n \eqdef \ppar*{Q_L^n , Q_1^n \ldots Q_J^n , Q_R^n}$ and the Euler time-step of this solution by an interval $\Delta t^n$ by $\mathscr{E}\ppar*{\pbrc{Q}^n}$. The Euler steps have a maximal time-step size given by the CFL condition for the bulk scheme, \eg for a finite volume scheme
\begin{equation}\label{eqn:FV_CFL}
	0 \leq \Delta t^n \max_j \ppar*{\frac{\max\abs{\lambda_j^{n} - \dot{x}_j^n}}{\Delta x_j^n}} \leq C_{\max}
\end{equation}
where $\max\abs{\lambda_j^{n} - \dot{x}_j^n}$ is the maximal characteristic speed inside of cell $j$ at time $t^{n}$ relative to the speed of mesh deformation $\dot{x}_j$, and $C_{\max}$ is the maximal permitted Courant number. The time-step from $t^n$ to $t^{n+1}$ is composed of $S$ sub-steps, each sub-step taking the form
\begin{subequations}\label{eqn:FV_RK}\begin{align}
	\pbrc*{Q}^{n [s+1]} &= \eta^{[s]} \pbrc*{Q}^{n [0]} + \ppar*{1-\eta^{[s]}} \cdot \mathscr{E}\ppar*{\pbrc*{Q}^{n[s]}}
	&
	\text{at}
	&&
	t^{n [s+1]} &= \eta^{[s]} t^{n [0]} + \ppar*{1-\eta^{[s]}} \cdot \ppar*{t^{n[s]} + \Delta t^n},
\end{align}\end{subequations}
with the intermediate values $\pbrc*{Q}^{n [s]}$ approximating the solution at time $t^{n[s]}$, $\pbrc*{Q}^{n [0]} \eqdef \pbrc*{Q}^{n}$, and $\pbrc*{Q}^{n+1} \eqdef \pbrc*{Q}^{n [S]}$. The values of $S$ and $\eta^{[s]}$ for time-stepping at different orders are given in \cref{tab:RKcoeff}. In order to evolve both the bulk and boundary values, we must establish a suitable Euler-like time-step for the algebraic boundary conditions (\cref{sec:TDP}).

Above we remarked that the scheme \eqref{eqn:FV_RK} is popular because it permits extrapolation of properties from Euler's method to higher order methods. We now make this more concrete.
\begin{theorem}	\label{thm:FV_time_step_stability_full}
	Let $z^{n[s]} \in Z$ be the consequence of time-stepping $z^{n[0]}$ using \eqref{eqn:FV_RK}, where $Z$ is a real vector space with a semi-norm $\norm{\omitdummy}$, and let $\hat{z}\ppar*{t}$ be some differentiable function from $t \in [t^{n},t_{\star}^n]$ to $Z$ where $t_{\star}^n = \max_{s<S}\ppar*{t^{n[s]} + \Delta t^n}$ with $\dv*{\hat{z}}{t}$ possessing a Lipschitz constant $L_z$. If, for some constant $\tilde{K} \geq 0$,
	\begin{equation}
		\norm*{\mathscr{E}\ppar*{z^{n[s]}} - \hat{z}\ppar*{t^{n[s]} + \Delta t^n}} \leq \norm*{z^{n[s]} - \hat{z}\ppar*{t^{n[s]}}} + \tilde{K} 
	\end{equation}
	for all $s \in \{0,1,\ldots,S-1\}$, then
	\begin{equation}
		\norm*{z^{n+1} - \hat{z}\ppar*{t^{n+1}}} \leq \norm*{z^{n} - \hat{z}\ppar*{t^{n} \vphantom{t^{n+1}}}} + \tilde{K}  + \frac{7+9S}{12+4S} L_z \ppar*{\Delta t^n}^2.
	\end{equation}
\end{theorem}
A proof is provided in \cref{sec:RK_contraction_proof}. As a brief aside, we observe that \cref{thm:FV_time_step_stability_full} allows us to confirm the TVD property of \eqref{eqn:FV_RK} by taking $z^n = \pbrc*{Q}^n$, $\hat{z} = 0$, and $\norm{\omitdummy} = \TV\ppar*{\omitdummy}$. Thus $\TV\ppar*{\mathscr{E}\ppar*{\pbrc*{Q}}} \leq \TV\ppar*{\pbrc*{Q}}$ implies  $\TV(\pbrc*{Q}^{n+1}) \leq \TV(\pbrc*{Q}^{n})$ as was stated in \cite{ar_Shu_1988} (see their proposition 2.1 and associated remark).

%% file: Sections/Spatial.tex
\section{Implementing and discretising extrapolation}	\label{sec:spatial}

\input{Sections/Extrapolation}
\input{Sections/FiniteDifferences}
\input{Sections/Diffusion}

%% file: Sections/Extrapolation.tex
\subsection{Extrapolation and forcing}	\label{sec:Extrapolation}

To evolve the boundary values $Q_R$ and $\dot{x}_R$ we combine extrapolation from the bulk with the characteristic equations. We modify the $m^\nth$ characteristic equation by introducing a forcing term, so that at the boundary
\begin{align}
	\label{eqn:hyp_extp_characteristic}
	l^{(m)} \dv{Q_R}{t} &= \ppar*{\dot{x}_R-\lambda^{(m)}} \cdot l^{(m)} Q_{xR} + l^{(m)} \cdot \Psi' + D^{(m)} l^{(m)} \cdot \ppar*{\hat{Q}_R - Q_R}.
\end{align}
for some $D^{(m)} \geq 0$. The value of $\hat{Q}_R$ is constructed by extrapolation, and the new term forces the value of $Q_R$ to tend towards $\hat{Q}_R$ over time, thus the simulated solution will satisfy \cref{eqn:boundary_limit}. The weak enforcement of $Q_R \rightarrow \hat{Q}_R$ will improve accuracy when the extrapolation is more accurate than the evolution, such as for static characteristics. The structure of the forcing term means that in the linear ($F = AQ$, $A$ constant) homogeneous ($\Psi = 0$) case, each of the characteristic fields $l^{(m)}Q$ evolves independently.

We now combine our modified characteristic equation with the boundary conditions. The system to solve for the boundary values is, for DA conditions,
\begin{align}\label{eqn:BC_extp_system_general}
	\pbrk*{\begin{array}{@{}cc@{}}
		L^+				&	0			\\
		B_Q				&	B_x
	\end{array}}
	\dv{}{t}
	\begin{bmatrix}
		Q_R	\\
		\dot{x}_R
	\end{bmatrix}
	&=
	\begin{bmatrix}
		\ppar*{\dot{x}_R I-\Lambda^+ } L^+ Q_{xR} + L^+ \Psi'+ D^+ L^+ \ppar*{\hat{Q}_R - Q_R}	\\
		b
	\end{bmatrix},	&
	g &= 0,
\end{align}
where $D^+ \eqdef \diag \pbrk*{ D^{(M_i+1)}, \ldots,  D^{(M)}}$.

We also consider the case of non-reflecting conditions. While in practice they may be imposed by specialist means, including them is informative and a check of robustness. It was proven in \cite{ar_Hedstrom_1978} that, in regions where the solution is continuous and the source term is zero, it is appropriate to enforce $l^{(m)} \pdv*{Q}{x} = 0$ as a boundary condition, which is also appropriate for a system with non-zero source terms. Thus the system for the boundary values is
\begin{align}\label{eqn:BC_extp_system_reflect}
	\pbrk*{\begin{array}{@{}cc@{}}
		L^+				&	0			\\
		L^-				&	0			\\
		B_Q				&	B_x
	\end{array}}
	\dv{}{t}
	\begin{bmatrix}
		Q_R	\\
		\dot{x}_R
	\end{bmatrix}
	&=
	\begin{bmatrix}
		\ppar*{\dot{x}_R I-\Lambda^+ } L^+ Q_{xR} + L^+ \Psi'+ D^+ L^+ \ppar*{\hat{Q}_R - Q_R}	\\
		L^- \Psi' + D^- L^- \ppar*{\hat{Q}_R - Q_R}	\\
		b
	\end{bmatrix},	&
	g &= 0,
\end{align}
where $D^- \eqdef \diag \pbrk*{ D^{(1)}, \ldots,  D^{(M_i)}}$, and we include a single DA condition to specify the boundary location, speed, or acceleration


To discretise the system we require expressions for $Q_{xR}$ and $\hat{Q}_R$. For the purpose of proving results we use general $k+1$ point expressions which can achieve up to $(k+1)^\nth$ order error. For extrapolation we take
\begin{align}
	\hat{Q}_R &= \sum_{j = J-k}^{J} \kappa_j Q_j = \kappa_J Q_J + \ldots + \kappa_{J-k} Q_{J-k},
\end{align}
and define a gradient function
\begin{align}
	\Gamma\ppar*{\bar{Q}_R,\bar{Q}_J \ldots \bar{Q}_{J-k+1}} &\eqdef  \sum_{j = J-k+1}^{J+1} \frac{\gamma_j \bar{Q}_j}{\Delta x_{J+1/2}}
	= \frac{1}{\Delta x_{J+1/2}} \ppar*{ \gamma_R \bar{Q}_R + \ldots + \gamma_{J-k+1} \bar{Q}_{J-k+1} },
\end{align}
($\gamma_{J+1} \equiv \gamma_R$, $\bar{Q}_{J+1} \equiv \bar{Q}_R$, $\bar{Q}_{\omitdummy}$ are dummy variables)
from which we define the following gradients
\begin{align}
	Q_{xR} &= \Gamma\ppar*{Q_R,Q_J \ldots Q_{J-k+1}},
	&
	\hat{Q}_{xR} &= \Gamma\ppar*{\hat{Q}_R,Q_J \ldots Q_{J-k+1}}.
\end{align}
The expression $Q_{xR}$ is to be used in the scheme, $\hat{Q}_{xR}$ will be used in our analysis. We next present extrapolations and derivatives to be used at different orders.

\paragraph{Zeroth order extrapolation}

For $k=0$ we have
\begin{align} \label{eqn:BCSP_extrap0}
	\hat{Q}_R &= Q_{J}, &
	\Gamma\ppar*{\bar{Q}_R} &= 0 = Q_{xR}.
\end{align}
Introduces first order error in smooth regions, but does not cause problems around shocks and gradient discontinuities.

\paragraph{First order extrapolation}

\begin{figure}[tp!]
	\centering
	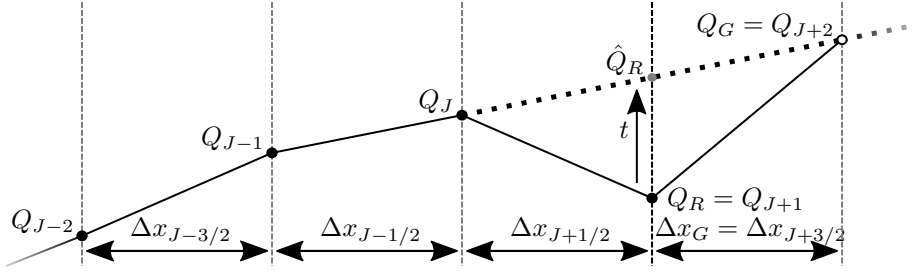
	\caption{Construction of the derivative and extrapolation at the right domain end. The dotted line indicates the gradient that we limit towards when using first order extrapolation, and the grey circle is the value of $\hat{Q}_R$, which $Q_R$ will limit towards as $t \rightarrow +\infty$ for a static field provided the bulk values do not change. Also plotted is the ghost point with value $Q_G$ used in the construction of the second derivative in \cref{sec:Diffusion}.}
	\label{fig:Derivative_Discrete}
\end{figure}

For $k=1$, depicted in \cref{fig:Derivative_Discrete}, and employed in \cref{sec:TP}. We define
\begin{align}\label{eqn:findif_for_scheme}
	\pbrk*{ \pdv{Q}{x} }_{j+1/2} &\eqdef \frac{Q_{j+1}-Q_{j}}{\Delta x_{j+1/2}},
\end{align}
the nearest neighbour finite difference derivatives. The extrapolation and derivative for the boundary conditions are
\begin{align} \label{eqn:BCSP_extrap1}
	\hat{Q}_R &= Q_{J} + \Delta x_{J+1/2} \pbrk*{ \pdv{Q}{x} }_{\mathrlap{J-1/2}},	&
	\Gamma\ppar*{\bar{Q}_R,\bar{Q}_J} &= \frac{\bar{Q}_R - \bar{Q}_J}{\Delta x_{J+1/2}},	&
	Q_{xR} &= \pbrk*{ \pdv{Q}{x} }_{\mathrlap{J+1/2}}.
\end{align}
Introduces second order error in smooth regions. When a gradient discontinuity collides with the domain end it introduces first order error. As a shock collides with the domain end zeroth order error will be generated, i.e. error which does not improve with resolution, but because it is localised the error advected into the bulk will be $\order{\Delta x}$. These errors are the same with and without forcing, and are a consequence of the polynomial interpolation used to deduce $Q_{xR}$ (and $\hat{Q}_R$) being less accurate in these cases.

\paragraph{Higher order extrapolation} For $k \geq 2$ there is a distinction between finite difference and finite volume schemes: for finite difference schemes the interpolation should take the value $Q_j$ at the point $x_j$; while for finite volume schemes it should integrate to $Q_j$ over the $j^\nth$ cell. It may be desirable to adjust the order of extrapolation depending on the local smoothness in a similar manner to weighted essentially non-oscillatory (WENO) schemes.

To close we overview our algorithm so far. First we compute the value of $\hat{Q}_R$ by extrapolation from bulk and we calculate the derivative, in \cref{sec:TP} we use \cref{eqn:BCSP_extrap1}. We also require an expression for $D^{(m)}$, and we discuss suitable values in \cref{sec:bc_discrete}. To assemble the system \cref{eqn:BC_extp_system_general} or \cref{eqn:BC_extp_system_reflect} for non-linear systems we additionally need evaluations of the eigenvalues and eigenvectors local to $x=x_R$, discussed in \cref{tp:TP_critwall}. To evolve we require an Euler time-step able to handle algebraic constraints, discussed in \cref{sec:TDP}. 

%% file: Figures/BoundaryDerivatives.pdf_tex
\begingroup%
  \makeatletter%
  \providecommand\color[2][]{%
    \errmessage{(Inkscape) Color is used for the text in Inkscape, but the package 'color.sty' is not loaded}%
    \renewcommand\color[2][]{}%
  }%
  \providecommand\transparent[1]{%
    \errmessage{(Inkscape) Transparency is used (non-zero) for the text in Inkscape, but the package 'transparent.sty' is not loaded}%
    \renewcommand\transparent[1]{}%
  }%
  \providecommand\rotatebox[2]{#2}%
  \newcommand*\fsize{\dimexpr\f@size pt\relax}%
  \newcommand*\lineheight[1]{\fontsize{\fsize}{#1\fsize}\selectfont}%
  \ifx\svgwidth\undefined%
    \setlength{\unitlength}{340.15748031bp}%
    \ifx\svgscale\undefined%
      \relax%
    \else%
      \setlength{\unitlength}{\unitlength * \real{\svgscale}}%
    \fi%
  \else%
    \setlength{\unitlength}{\svgwidth}%
  \fi%
  \global\let\svgwidth\undefined%
  \global\let\svgscale\undefined%
  \makeatother%
  \begin{picture}(1,0.29166667)%
    \lineheight{1}%
    \setlength\tabcolsep{0pt}%
    \put(0,0){\includegraphics[width=\unitlength,page=1]{BoundaryDerivatives.pdf}}%
    \put(0.81666667,0.03333333){\color[rgb]{0,0,0}\makebox(0,0)[t]{\lineheight{1.25}\smash{\begin{tabular}[t]{c}$\Delta x_G = \Delta x_{J+3/2}$\end{tabular}}}}%
    \put(0,0){\includegraphics[width=\unitlength,page=2]{BoundaryDerivatives.pdf}}%
    \put(0.60833333,0.03333333){\color[rgb]{0,0,0}\makebox(0,0)[t]{\lineheight{1.25}\smash{\begin{tabular}[t]{c}$\Delta x_{J+1/2}$\end{tabular}}}}%
    \put(0.4,0.03333333){\color[rgb]{0,0,0}\makebox(0,0)[t]{\lineheight{1.25}\smash{\begin{tabular}[t]{c}$\Delta x_{J-1/2}$\end{tabular}}}}%
    \put(0.19166665,0.03333333){\color[rgb]{0,0,0}\makebox(0,0)[t]{\lineheight{1.25}\smash{\begin{tabular}[t]{c}$\Delta x_{J-3/2}$\end{tabular}}}}%
    \put(0,0){\includegraphics[width=\unitlength,page=3]{BoundaryDerivatives.pdf}}%
    \put(0.68333333,0.14166667){\color[rgb]{0,0,0}\makebox(0,0)[rt]{\lineheight{1.25}\smash{\begin{tabular}[t]{r}$t$\end{tabular}}}}%
    \put(0.725,0.06666667){\color[rgb]{0,0,0}\makebox(0,0)[lt]{\lineheight{1.25}\smash{\begin{tabular}[t]{l}$Q_R=Q_{J+1}$\end{tabular}}}}%
    \put(0.49166667,0.175){\color[rgb]{0,0,0}\makebox(0,0)[rt]{\lineheight{1.25}\smash{\begin{tabular}[t]{r}$Q_{J}$\end{tabular}}}}%
    \put(0.28333333,0.13333333){\color[rgb]{0,0,0}\makebox(0,0)[rt]{\lineheight{1.25}\smash{\begin{tabular}[t]{r}$Q_{J-1}$\end{tabular}}}}%
    \put(0.90833333,0.25833333){\color[rgb]{0,0,0}\makebox(0,0)[rt]{\lineheight{1.25}\smash{\begin{tabular}[t]{r}$Q_{G}=Q_{J+2}$\end{tabular}}}}%
    \put(0.07499999,0.04166667){\color[rgb]{0,0,0}\makebox(0,0)[rt]{\lineheight{1.25}\smash{\begin{tabular}[t]{r}$Q_{J-2}$\end{tabular}}}}%
    \put(0.7,0.21666667){\color[rgb]{0,0,0}\makebox(0,0)[rt]{\lineheight{1.25}\smash{\begin{tabular}[t]{r}$\hat{Q}_{R}$\end{tabular}}}}%
  \end{picture}%
\endgroup%

%% file: Sections/FiniteDifferences.tex
\subsection{Bounding the forcing} \label{sec:bc_discrete}

In this section we construct bounds on $D^{(m)}$ to ensure that the system is stable under the scheme \cref{eqn:FV_RK}. We do this for a linear flux so that $A$ is constant (the case of a non-linear flux is considered in \cref{sec:TP}), with a Lipschitz continuous source so that
\begin{align}
	\abs*{l^{(m)} \cdot \ppar*{\Psi'\ppar*{Q_1,x_R,t} - \Psi'\ppar*{Q_2,x_R,t} }}
	&\leq P^{(m)} \abs*{l^{(m)} \cdot \ppar*{Q_1 - Q_2 }}
	&\text{for all}&&
	Q_1,Q_2
\end{align}
where $P^{(m)}>0$ are constants. To simplify expressions, here and in the rest of this section we work at time $t = t^{n[s]}$, and we employ an Euler time-step of size $\Delta t = \Delta t^n$. We will assume that $Q_{J-k} \ldots Q_{J}$ converge as $\Delta x \rightarrow 0$; in particular, we treat the bulk values as independent of the boundary value $Q_R$ (\eg they come from the exact solution). Convergence in the fully coupled case where the bulk values are influenced by the boundary is demonstrated in \cref{sec:TP}. If the bulk scheme is stable and the solution is locally Lipschitz, the evolution of $\hat{Q}_R$ will satisfy the assumptions of $\hat{z}$ in \cref{thm:FV_time_step_stability_full}. Without knowing the details of the bulk scheme we cannot establish a value for $\mathscr{E}\ppar*{\hat{Q}_R}$. However, given that the bulk scheme is consistent (\ie vanishing truncation error) we can evolve according to the bulk values, which means that for outgoing characteristics we use the gradient based on extrapolation $\hat{Q}_{xR}$. For incoming characteristics, if they are to be evolved using the characteristic equations then this can only be because we impose a non-reflecting condition $l^{(m)} Q_{xR}=0$. Thus we approximate the time-stepping of $\hat{Q}_{R}$ by $\hat{Q}_{R}^{\text{step}}$, where
\begin{equation}	\label{eqn:BCSP_centered_target_evolve}
	l^{(m)} \hat{Q}_{R}^{\text{step}} \eqdef l^{(m)}\hat{Q}_{R} + \Delta t \, l^{(m)} \ppar*{ -\max \ppar*{\lambda^{(m)} - \dot{x}_R,0} \cdot \hat{Q}_{xR} + \Psi'\ppar*{\hat{Q}_{R},x_R,t}}.
\end{equation}
Because the bulk scheme is consistent and at least first order we can say that
\begin{equation}	\label{eqn:BCSP_centered_target_difference}
	\norm*{\hat{Q}_{R}^{\text{step}} - \mathscr{E}\ppar*{\hat{Q}_{R}}} = \order{\ppar*{\Delta t}^2,\ppar*{\Delta t}^{k+1}}.
\end{equation}
The boundary value under Euler time-stepping of \cref{eqn:hyp_diff_characteristic} is
\begin{equation}	\label{eqn:BCSP_value_evolve}
	\mathscr{E}\ppar*{ l^{(m)} Q_{R} } = l^{(m)} Q_{R} + \Delta t \, l^{(m)} \cdot \ppar*{ - \max \ppar*{\lambda^{(m)} - \dot{x}_R ,0} \cdot Q_{xR} + \Psi'\ppar*{Q_{R},x_R,t} + D^{(m)} \cdot \ppar*{\hat{Q}_R - Q_R} }.
\end{equation}

\begin{theorem}\begin{subequations} \label{eqn:BCSP_centered_convergence_single}
	If $m \in \{1,2,\ldots,M\}$, \cref{eqn:BCSP_centered_target_difference} is satisfied, and
	\begin{gather}	\label{eqn:BCSP_centered_convergence_single_condition}
		P^{(m)} - \max \ppar*{\lambda^{(m)} - \dot{x}_R ,0} \frac{\gamma_R}{\Delta x_{J+1/2}} \leq D^{(m)} \leq \frac{2}{\Delta t} - \max \ppar*{\lambda^{(m)} - \dot{x}_R ,0} \frac{\gamma_R}{\Delta x_{J+1/2}} - P^{(m)}
	\shortintertext{then}
		\label{eqn:BCSP_centered_convergence_single_result}
		\abs*{\mathscr{E} \ppar*{l^{(m)}Q_{R}} - \mathscr{E}\ppar*{l^{(m)}\hat{Q}_{R}}} \leq \abs*{l^{(m)}\ppar*{Q_{R} - \hat{Q}_{R}}} + \order{\ppar*{\Delta t}^2,\ppar*{\Delta t}^{k+1}}.
	\end{gather}
\end{subequations}\end{theorem}
\begin{proof}
	First we observe that
	\begin{align*}
		\hat{Q}_{xR} - Q_{xR} &= \gamma_R \frac{\hat{Q}_{R} - Q_{R}}{\Delta x_{J+1/2}}.
	\end{align*}
	Thus, denoting $\lambda = \max (\lambda^{(m)} - \dot{x}_R ,0)$,
	\begin{align*}
		\mathscr{E} \ppar*{l^{(m)}Q_{R}} - l^{(m)}\hat{Q}_{R}^{\text{ step}}
		&=
		\ppar*{ 1 - \frac{\lambda \gamma_R \Delta t}{\Delta x_{J+1/2}} - D^{(m)} \Delta t } \cdot l^{(m)} \cdot \ppar*{ Q_{R} - \hat{Q}_{R} }
		+ \Delta t \, l^{(m)} \cdot \ppar*{ \Psi'(Q_{R},x_R,t) - \Psi'\ppar*{\hat{Q}_{R},x_R,t} }
	\shortintertext{and}
		\abs*{\mathscr{E} \ppar*{l^{(m)}Q_{R}} - \mathscr{E}\ppar*{l^{(m)}\hat{Q}_{R}}}
		&\leq
		\ppar*{\abs*{ 1 - \frac{\lambda \gamma_R \Delta t}{\Delta x_{J+1/2}} - D^{(m)} \Delta t} + \Delta t P^{(m)}} \abs*{l^{(m)} \cdot \ppar*{ Q_{R} - \hat{Q}_{R} }} + \order{\ppar*{\Delta t}^2}.
	\end{align*}
	Therefore, by \cref{thm:FV_time_step_stability_full}, for \eqref{eqn:BCSP_centered_convergence_single_result} to be true it is sufficient that
	\begin{align*}
		\abs*{ D^{(m)} \Delta t +  \frac{\lambda \gamma_R \Delta t}{\Delta x_{J+1/2}} - 1 } + \Delta t P^{(m)} &\leq 1
	\end{align*}
	which is equivalent to \cref{eqn:BCSP_centered_convergence_single_condition}.
\end{proof}


The presence of a non-zero lower bound for $D^{(m)}$ is perhaps surprising. We see that it depends on the Lipschitz bounds for $\Psi'$, denoted $P^{(m)}$, revealing that this lower bound comes from variation of the source strength with $Q$. Indeed, if $Q_J,Q_{J-1}\ldots$ are increasing, $Q_R$ is smaller than $\hat{Q}_R,Q_J,Q_{J-1},\ldots$, and the source strength increases in $Q$, then $\hat{Q}_R$ will grow faster than $Q_R$ causing $\hat{Q}_R-Q_R$ to grow, and to close the gap we require $D^{(m)} = \order{1}$. Similarly, the upper bound can be understood as a result of advective, forcing, and diffusive effects acting on the difference between $Q_R$ and $\hat{Q}_R$, and so some compromise between these effects is required.

As $\Delta x \rightarrow 0$, the bounds \cref{eqn:BCSP_centered_convergence_single_condition} can be satisfied by requiring that
\begin{align}	\label{eqn:BCSP_onesided_convergence_weakened}
	0 &\leq D^{(m)} < \frac{1}{\Delta x_{J+1/2}} \ppar*{ \frac{1}{C_{\max}} - \gamma_R } \max_{m'}\abs*{\lambda^{(m')} - \dot{x}_R}
	\qquad \text{for all } m,
\end{align}
($\gamma_R=1$ for first order extrapolation) where we used that the finite volume CFL condition \cref{eqn:FV_CFL} implies
\begin{align}
	\Delta t  \leq C_{\max} \frac{\Delta x_J}{\max_{m'}\abs*{\lambda^{(m')} - \dot{x}_R}} = C_{\max} \frac{2 \Delta x_{J+1/2}}{\max_{m'}\abs*{\lambda^{(m')} - \dot{x}_R}}.
\end{align}
In \cref{eqn:BCSP_onesided_convergence_weakened} the upper bound is strict so that \eqref{eqn:BCSP_centered_convergence_single_condition} is satisfied at a finite resolution; the lower bound is non-strict to include the forcing free case. When computing the maximal time-step we have assumed that it is limited by the dynamics at the boundary. In a non-linear system this may not be the case and a larger amount of forcing may be stable. Thus, \eqref{eqn:BCSP_onesided_convergence_weakened} should be sufficient for stability for most problems, as demonstrated in \cref{tp:TP_alphawave}.

Provided that no catastrophic interaction with the method for the bulk occurs, our results are sufficient to guarantee convergence for regions in which the high resolution limit yields locally constant gradient, eigenvalues, and eigenvectors. The regions where this will not happen are local to shocks, and close to points where the system becomes degenerate, which we investigate in \cref{sec:TP}.

%% file: Sections/Diffusion.tex
\subsection{Reinterpreting extrapolation as diffusion}	\label{sec:Diffusion}

We now take a brief aside and discuss an interpretation of the correction term introduced in \cref{eqn:hyp_extp_characteristic} in terms of diffusion. For this purpose we introduce a ghost point at location $x=x_G$, $\Delta x_G \eqdef x_G - x_R$, with value $Q_G$ (see \cref{fig:Derivative_Discrete}). We construct finite difference second derivatives
\begin{align}
	\pbrk*{ \pdv[2]{Q}{x} }_{j} &\eqdef \frac{2}{\Delta x_{j+1/2} + \Delta x_{j-1/2}} \ppar*{ \pbrk*{ \pdv{Q}{x} }_{j+1/2} - \pbrk*{ \pdv{Q}{x} }_{j-1/2} }
\end{align}
where, for simplicity of notation, we take $Q_{J+2} = Q_{G}$, and $\Delta x_{J+3/2} = \Delta x_{G}$. Choosing
\begin{align}
	\label{eqn:hyp_diff_second}
	Q_G &= \hat{Q}_R + \Delta x_G \frac{\hat{Q}_R - Q_J}{\Delta x_{J+1/2}}
	&\text{we obtain}&&
	\pbrk*{ \pdv[2]{Q}{x} }_{J+1} &= \frac{2}{\Delta x_{J+1/2}\Delta x_{G}} \ppar*{ \hat{Q}_R - Q_{R} },
\end{align}
the value of $Q_G$ lying on the linear extrapolation from $Q_J$ through $\hat{Q}_R$ as depicted in \cref{fig:Derivative_Discrete}. Thus we may consider the approach developed in \cref{sec:Extrapolation,sec:bc_discrete} as a careful discretisation of the equation
\begin{align}
	\label{eqn:hyp_diff_characteristic}
	l^{(m)}\dv{Q_R}{t} &= \ppar*{\dot{x}_R-\lambda^{(m)}} l^{(m)} Q_{xR} + l^{(m)} \Psi' + \tilde{D}^{(m)} l^{(m)} \pdv[2]{Q}{x}
	&\text{where}&&
	\tilde{D}^{(m)} &= \frac{\Delta x_{J+1/2}\Delta x_{G}}{2} D^{(m)},
\end{align}
which is the projection onto the boundary of the system
\begin{align}
	\label{eqn:hyp_diff_system}
	\pdv{Q}{t} + \pdv{}{x}\ppar*{F} &= \Psi + R \tilde{D} L \pdv[2]{Q}{x}	
	&\text{where}&&
	\tilde{D} &= \diag \pbrk*{\tilde{D}^{(1)} , \ldots , \tilde{D}^{(M)}}.
\end{align}
As $\Delta x \rightarrow 0$, $\tilde{D}^{(m)} \rightarrow 0$, thus our approach is consistent with the vanishing viscosity limit in the bulk with an unusual diffusivity matrix. See \cite{ar_Bianchini_2005,ar_Christoforu_2006} for results on the vanishing viscosity solution with diffusivity proportional to the identity, and \cite{ar_Bianchini_2003,ar_Bianchini_2016} for Riemann and boundary Riemann problems with a broader class of diffusivity matrices.

To finish our brief discussion, we present the expressions for the second derivative produced by the two methods of extrapolation considered in \cref{sec:bc_discrete}. The expressions are unusual, but we can be confident that they do represent a second derivative by their construction. Zeroth order \cref{eqn:BCSP_extrap0} and first order \cref{eqn:BCSP_extrap1} extrapolation yield, respectively, 
\begin{align}
	\Delta x_{G} \pbrk*{ \pdv[2]{Q}{x} }_{J+1}  &= - 2 \pbrk*{ \pdv{Q}{x} }_{\mathrlap{J+1/2}} , 
	&\text{and}&&
	\Delta x_{G} \pbrk*{ \pdv[2]{Q}{x} }_{J+1}  &= - \ppar*{ \Delta x_{J+1/2} + \Delta x_{J-1/2} }\pbrk*{ \pdv[2]{Q}{x} }_{\mathrlap{J}}.
\end{align}

%% file: Sections/Temporal.tex
\section{Temporal discretization and projection}	\label{sec:TDP}

\input{Sections/Temporal_Intro}
\input{Sections/Temporal_Manifold}
\input{Sections/Temporal_Projection}
\input{Sections/Temporal_RKNR}
\input{Sections/Temporal_Comparison}

%% file: Sections/Temporal_Intro.tex
We want to simulate the system for  $v : t \mapsto \mathbb{R}^{\tilde{M}}$
\sidebysidesubequations{eqn:TDP_DAE_full}{
	\tilde{B}(v,t) \cdot \dv{v}{t} = \tilde{b}(v,t),
}{eqn:TDP_DAE_diff}{}{
	\tilde{g}(v,t) = 0,
}{eqn:TDP_DAE_alge}
where $\tilde{B} : (v,t) \mapsto \mathbb{M}\ppar*{\tilde{M}-\tilde{M}_A,\tilde{M}}$, $\tilde{b} : (v,t) \mapsto \mathbb{R}^{\tilde{M}-\tilde{M}_A}$, and $\tilde{g} : (v,t) \mapsto \mathbb{R}^{\tilde{M}_A}$. The system \cref{eqn:TDP_DAE_full} contains the full set of spatially discrete equations, including both left and right boundary conditions and the full set of evolution equations for the bulk values. However, in the schemes presented in this section, all values other than $\ppar{Q_R^T , \dot{x}_R}^T$ (and the equivalent at $x_L$) are evolved using the RK scheme \cref{eqn:FV_RK} and it is only the evolution of the boundary values that is modified. For this reason we may alternatively interpret the system as being for $v = \ppar{Q_R^T , \dot{x}_R}^T$ with all other values known functions of $t$, in which case $\tilde{M} = M+1$ and $\tilde{M}_A = M_A$.

It is the presence of algebraic constraints that complicates the evolution of \eqref{eqn:TDP_DAE_full}; \eqref{eqn:TDP_DAE_alge} defines a manifold on which $v$ evolves according to \eqref{eqn:TDP_DAE_diff}. A variety of numerical methods exist for differential equations on manifolds; we discuss those which are compatible with the RK scheme \eqref{eqn:FV_RK}.

%% file: Sections/Temporal_Manifold.tex
\subsection{Evolution in a local coordinate system}	\label{sec:TDP_local}

This scheme operates on the manifold, using a local coordinate system $(w,t)$ where $w\in\mathbb{R}^{\tilde{M} - \tilde{M}_a}$ and the points on the manifold transform as $v = V(w,t)$, $w = W(v,t)$. Substituting this change of variables into \eqref{eqn:TDP_DAE_diff} yields
\begin{subequations} \label{eqn:TDP_local}
\begin{gather}
	\tilde{B}_w (w,t) \cdot \dv{w}{t} = \tilde{b}_w (w,t)
\shortintertext{where}
	\begin{aligned}
		\tilde{B}_w (w,t) &= \tilde{B} \ppar*{(V(w,t),t} \cdot \pdv{V}{w} (w,t),
		&\text{and}&&
		\tilde{b}_w (w,t) &= \tilde{b} \ppar*{V(w,t),t} - \tilde{B}\ppar*{V(w,t),t} \cdot \pdv{V}{t} (w,t).
\end{aligned}
\end{gather}\end{subequations}
We can use this change of coordinates to perform each sub-step of the RK scheme in the local coordinate system. That is to say, we transform from the value at any given sub-step $v^{n[s]}$ to $w^{n[s]} = W\ppar{v^{n[s]},t^{n[s]}}$, perform the sub-step on \eqref{eqn:TDP_local} using \eqref{eqn:FV_RK} to obtain $w^{n[s+1]}$, and then transform back as $v^{n[s+1]} = V\ppar{w^{n[s+1]},t^{n[s+1]}}$. The full RK step must all be performed in the same sufficiently smooth local coordinate system in order to produce the order of accuracy claimed.

This scheme is the most accurate of those we will present here, perfectly resolving the manifold. However, it requires an atlas of local coordinate systems, each covering enough of the manifold to take any permitted time-step. This is a major drawback, because a suitable atlas is usually very hard to produce, rendering this approach impractical. If an atlas can be found then this approach is recommended. See, for example, \cite{bk_Hairer_SDAP} for more information on this scheme, which we term RKLC (Runge-Kutta in local coordinates).

%% file: Sections/Temporal_Projection.tex
\subsection{Runge-Kutta with projection}	\label{sec:TDP_projection}

If we lack an atlas then we operate in $v$ instead. A simple approach (used in \cite{ar_Thompson_1987,ar_Thompson_1990,ar_Kim_2000,ar_Kim_2004}) is to differentiate \eqref{eqn:TDP_DAE_alge} with respect to $t$ so that the system becomes 
\begin{align}\label{eqn:TDP_DAE_ODEversion}
	\tilde{B}(v,t) \cdot \dv{v}{t} &= \tilde{b}(v,t),
	&
	\tilde{G}(v,t) \cdot \dv{v}{t} &= -\pdv{\tilde{g}}{t}(v,t),
\end{align}
where $\tilde{G} \eqdef \pdv*{\tilde{g}}{v}$, which can then be simulated using \eqref{eqn:FV_RK}. We term this scheme RK0 (Runge-Kutta without projection), and it suffers from drift-off error as discussed in \cref{sec:Intro}. For a hyperbolic system, the error generation local to certain times (such as at initiation, change of boundary condition, and interaction with a shock) can be $\order{1}$, causing the simulation to fail to converge on the correct boundary condition. One remedy is to use projection onto the manifold after each time-step (\eg \cite{bk_Ascher_CMODEDAE}), where by projection we mean a Newton-Raphson style iteration that brings us back onto the manifold. There are a variety of projection schemes that can be employed, for example
\begin{equation}
	v_{i+1}^{n} = v_{i}^{n} -
	\tilde{G}\ppar*{v_{i}^{n},t^{n}}^T \cdot \pbrk*{ \tilde{G}\ppar*{v_{i}^{n},t^{n}} \cdot \tilde{G}\ppar*{v_{i}^{n},t^{n}}^T }^{-1} \cdot \tilde{g}\ppar*{v_{i}^{n},t^{n}}
\end{equation}
where $v_{0}^{n}$ is the value produced by the RK scheme \eqref{eqn:FV_RK}. This iteration, at each step, projects towards the closet point at which $g=0$ as determined by the local gradients. If an $S^\text{th}$ order method is used and the functions in \eqref{eqn:TDP_DAE_ODEversion} are $S$ times differentiable, then the distance of $v_0^n$ from the manifold will be $\order{\ppar{\Delta t^{n-1}}^{S+1}}$. Indeed, since this is a Newton-Raphson (NR) style root-finding method we expect that $\tilde{g}\ppar*{v_{i}^{n},t^{n}} = \order{\ppar{\Delta t^{n-1}}^{(S+1) \cdot 2^i}}$, so long as $\tilde{G}$ has linearly independent rows. We term this scheme RKP1 (Runge-Kutta with post-step projection method 1).

An alternative projection scheme can be developed by examining the system \eqref{eqn:TDP_DAE_full} and thinking of the projection as a very small Euler time-step, after which we will be back on the manifold, that is
\begin{align}
	\tilde{B}(v,t) \cdot \Delta v &= \tilde{b}(v,t) \cdot \Delta t,
	&
	\tilde{g}(v + \Delta v,t+\Delta t) &=0.
\intertext{Limiting $\Delta t \rightarrow 0$, this yields the system}
	\tilde{B}(v,t) \cdot \Delta v &= 0,
	&
	\tilde{g}(v + \Delta v,t) &=0,
\end{align}
and solving using the Newton-Raphson method yields the iteration
\begin{equation}
	v_{i+1}^{n} = v_{i}^{n} -
	\pbrk*{\begin{array}{c}
		\tilde{B}\ppar*{v_i^{n},t^{n}}	\\
		\tilde{G}\ppar*{v_i^{n},t^{n}}
	\end{array}}^{-1}
	\pbrk*{\begin{array}{c}
		0	\\
		\tilde{g}\ppar*{v_i^{n},t^{n}}
	\end{array}}
\end{equation}
where we adjust the value of $B$ for each step to approximate tracking along a curve with $\tilde{B}(v,t) \cdot \dd v = 0$. Again, we expect that $\tilde{g} \ppar*{v_{i}^{n},t^{n}} = \order{\ppar{\Delta t^{n-1}}^{(S+1) \cdot 2^i}}$, so long as $\tilde{G}$ and $\tilde{B}$ have linearly independent rows. We term this scheme RKP2 (Runge-Kutta with post-step projection method 2).

It is also possible to use internal projection \cite{ln_Hairer_SDEM}, where we take the value after a sub-step as the initial value of the projection iteration, $v_0^{n[s]} \eqdef v^{n[s]}$, project $i$ times to yield $v_i^{n[s]}$, and then perform the Euler time-step as
\begin{equation}
	\mathscr{E}\ppar*{v^{n[s]}} = v^{n[s]} + \Delta t^n \eval*{\dv{v}{t}}_{\mathrlap{v = v_i^{n[s]}}}.
\end{equation}
This means that we time-step using a more accurate evaluation of $\dv*{v}{t}$, but note that the step starts from a location not necessarily on the manifold. After the time-step we use projection as in RKP1 and RKP2. We call these methods RKI1 and RKI2 (Runge-Kutta with internal projection method 1 and 2).

%% file: Sections/Temporal_RKNR.tex
\subsection{Runge-Kutta Newton-Raphson method} \label{sec:TDP_RKNR}

We can incorporate the ideas used to generate RKI2 into the numerical time-stepping algorithm itself, to construct an implicit method that uses iteration to solve $\tilde{g}=0$ in each time-step, while evolving the differential equations using Euler's method. This scheme was originally presented in \cite{ar_Skevington_F001_Draining}, and we reproduce it here in greater detail.

We wish to use the Runge-Kutta scheme \eqref{eqn:FV_RK}, and in the case $\tilde{M}_A=0$ the expression $\mathscr{E} \ppar*{v^{n[s]}}$ denotes the result of Euler time-stepping $v$ from time $t^{n[s]}$ to time $t^{n[s]}+\Delta t^n$. For our implicit scheme, we use the notation $\mathscr{E} \ppar*{v,\Delta v}$ for the approximate Euler time-step given that the previous attempt at the step yielded $v + \Delta v$. An iteration can then be defined by $\Delta v_{i+1}^{n[s]} \eqdef \mathscr{E} \ppar*{v^{n[s]},\Delta v_i^{n[s]}} - v^{n[s]}$, and we take $\mathscr{E}$ to be the limit denoted by $\mathscr{E}_{\infty}$, \ie
\begin{equation}
	\mathscr{E}_{\infty} \ppar*{v^{n[s]}} \eqdef v^{n[s]} + \lim_{i \rightarrow \infty} \Delta v_i^{n[s]}.
\end{equation}
The intermediate values for $v^{n[s+1]}$ are denoted by
\begin{equation}
	v_{i+1}^{n[s+1]} \eqdef \eta^{[s]} v^{n[0]} + \ppar*{1-\eta^{[s]}} \cdot \mathscr{E}\ppar*{v^{n[s]},\Delta v_{i}^{n[s]}} = \eta^{[s]} v^{n[0]} + \ppar*{1-\eta^{[s]}} \cdot \ppar*{v^{n[s]} + \Delta v_{i+1}^{n[s]}}
\end{equation}
thus $v^{n[s+1]} = \lim_{i \rightarrow \infty} v_i^{n[s+1]}$. We now construct the iteration, generating $\Delta v_{i+1}^{n[s]}$ from $\Delta v_{i}^{n[s]}$. The differential equation \eqref{eqn:TDP_DAE_diff} is time-stepped using Euler's method
\begin{subequations}
\begin{equation}
	\tilde{B} \ppar*{v^{n[s]},t^{n[s]}}  \cdot \Delta v_{i+1}^{n[s]} = \tilde{b} \ppar*{v^{n[s]},t^{n[s]}} \cdot \Delta t^n.
\end{equation}
Iteration is solely used to find the value of $\Delta v^{n[s]}$ which solves the equation $\tilde{g}\ppar{(v^{n[s+1]},t^{n[s+1]}} = 0$, so that we remain on the manifold $\tilde{g}=0$ after the sub-step. This is done using the Newton-Raphson method
\begin{equation}
	\ppar*{1-\eta^{[s]}} \cdot \tilde{G}\ppar*{v_i^{n[s+1]},t^{n[s+1]}} \cdot \ppar*{\Delta v_{i+1}^{n[s]} - \Delta v_i^{[s]}} = -\tilde{g}(v_i^{n[s+1]},t^{n[s+1]}).
\end{equation}
These can be combined into a single expression
\begin{equation}	\label{eq:BC_Euler_iteration}
	\mathscr{E}\ppar*{v^{n[s]},\Delta v_i^{n[s]}} = v^{n[s]} + 
	\pbrk*{ \begin{array}{c}
		\tilde{B}\ppar*{v^{n[s]} , t^{n[s]}}
		\\
		\ppar*{1-\eta^{[s]}} \cdot \tilde{G}\ppar*{v_i^{n[s+1]},t^{n[s+1]}}
	\end{array}}^{-1}
	\pbrk*{ \begin{array}{c}
		\tilde{b}\ppar*{v^{n[s]},t^{n[s]}} \cdot \Delta t
		\\
		\ppar*{1-\eta^{[s]}} \cdot \tilde{G}\ppar*{v_i^{n[s+1]},t^{n[s+1]}} \cdot \Delta v_i^{n[s]} - \tilde{g}\ppar*{v_i^{n[s+1]},t^{n[s+1]}}
	\end{array}}.
\end{equation}
\end{subequations}
We investigate the effectiveness of the method with respect to the differential and algebraic constraints separately. First, the differential constraints \eqref{eqn:TDP_DAE_diff} act in the directions spanned by the row-space of $\tilde{B}$, thus we should examine the projection of the error into this space. The projection is performed using the matrix
\begin{equation}
	\tilde{P} = \tilde{B}^T \ppar*{\tilde{B} \tilde{B}^T}^{-1} \tilde{B}.
\end{equation}
We define the truncation error of a single sub-step after $i$ iterations as
\begin{equation}
	\tau = \frac{1}{\Delta t^n}  \tilde{P}\ppar*{v^{n[s]},t^{n[s]}} \cdot \ppar*{ \mathscr{E}\ppar*{v^{n[s]},\Delta v_i^{n[s]}} - v\ppar*{t^{n[s]}+\Delta t^n} }
\end{equation}
denoting the exact solution at time $t$ by $v(t)$ with $v(t^{n[s]}) = v^{n[s]}$. We use the value of $\tilde{P}$ from the start of the time-step as Euler's method treats $\tilde{B}$ as constant over the time interval.

\begin{prop}\begin{subequations}
	The truncation error is
	\begin{align}\label{eqn:thm_RKNReuler_truncation}
		\tau = - \frac{\Delta t^n}{2} \tilde{B}^T \ppar*{\tilde{B} \tilde{B}^T}^{-1} \tilde{\alpha} + \order{\frac{1}{\Delta t^n}\norm*{v_i^{n[s+1]} - v}^3,\ppar*{\Delta t^n}^2 }
	\end{align}
	where all functions are evaluated at $v^{n[s]}$, $t^{n[s]}$, and
	\begin{align}
		\tilde{\alpha}_i &= \sum_j \pbrk*{ \pdv{\tilde{b}_i}{v_j}\dv{v_j}{t} - \pdv{\tilde{B}_{ij}}{t}\dv{v_j}{t} - \sum_k \pbrk*{ \pdv{\tilde{B}_{ij}}{v_k} \dv{v_j}{t} \dv{v_k}{t} } } + \pdv{\tilde{b}_i}{t}.
	\end{align}
\end{subequations}\end{prop}
\begin{proof}
	We first define
	\begin{align}
		\tilde{\beta}_i &= \sum_j \pbrk*{ 2\pdv{\tilde{g}_i}{t}{v_j} \dv{v_j}{t} + \sum_k \pbrk*{ \pdv{\tilde{g}_i}{v_j}{v_k} \dv{v_j}{t} \dv{v_k}{t} } } + \pdv[2]{\tilde{g}_i}{t},
	\end{align}
	The truncation error is found by Taylor expansion. We expand the exact solution as
	\begin{align}
		v\ppar*{t^{n[s]}+\Delta t^n} &= v + \Delta t^n \dv{v}{t} + \frac{\ppar*{\Delta t^n}^2}{2} \dv[2]{v}{t} + \order{\ppar*{\Delta t^n}^3}
	\end{align}
	where, by \cref{eqn:TDP_DAE_full}
	\begin{align}
		\tilde{B} \dv{v}{t} &= \tilde{b},
		&
		\tilde{B} \dv[2]{v}{t} &= \tilde{\alpha},
		&
		\tilde{G} \dv{v}{t} &= -\dv{\tilde{g}}{t},
		&\text{and}&&
		\tilde{G} \dv[2]{v}{t} &= -\tilde{\beta}.
	\end{align}
	Thus
	\begin{align}\begin{split}
		\tau ={}& \frac{\tilde{P}}{\Delta t^n}
		\left(\ppar*{
		v + 
		\pbrk*{ \begin{array}{c}
			\tilde{B}
			\\
			\ppar*{1-\eta^{[s]}} \cdot \tilde{G}\ppar*{v_i^{n[s+1]},t^{n[s+1]}}
		\end{array}}^{-1}
		\pbrk*{ \begin{array}{c}
			\tilde{b} \Delta t^n
			\\
			\ppar*{1-\eta^{[s]}} \cdot \tilde{G}\ppar*{v_i^{n[s+1]},t^{[s+1]}} \Delta v_i^{n[s]} - \tilde{g}\ppar*{v_i^{n[s+1]},t^{n[s+1]}}
		\end{array}}
		}
		\right. \\& \left.
		- \ppar*{
		v + \Delta t^n
		\pbrk*{ \begin{array}{c}\tilde{B}\\\tilde{G}\end{array}}^{-1}
		\pbrk*{ \begin{array}{c} \tilde{b}\\-\pdv*{\tilde{g}}{t}\end{array}}
		+\frac{(\Delta t^n)^2}{2}
		\pbrk*{ \begin{array}{c}\tilde{B}\\\tilde{G}\end{array}}^{-1}
		\pbrk*{ \begin{array}{c}\tilde{\alpha}\\-\tilde{\beta}\end{array}}
		}\right)
		+ \order{ \frac{1}{\Delta t^n}\norm*{v_i^{n[s+1]} - v}^3,\ppar*{\Delta t^n}^2 }.
	\end{split}\end{align}
	To obtain \cref{eqn:thm_RKNReuler_truncation} we use that
	\begin{equation}
		\tilde{P} 
		\pbrk*{ \begin{array}{c}\tilde{B}\\\tilde{Z}\end{array}}^{-1}
		=
		\tilde{B}^T \ppar*{\tilde{B} \tilde{B}^T}^{-1} \tilde{B}
		\pbrk*{ \begin{array}{c}\tilde{B}\\\tilde{Z}\end{array}}^{-1}
		=
		\pbrk*{ \begin{array}{cc}\tilde{B}^{T} \ppar*{\tilde{B} \tilde{B}^T}^{-1} & 0\end{array}}.
	\end{equation}
	for any matrix $\tilde{Z}$ so that the inverse exists.
\end{proof}

Thus the Euler's method portion of the scheme performs exactly as we would expect, and the truncation error for the ODE is $\order{\Delta t}$. We next investigate the convergence to the manifold.

\begin{prop}
	The value of the manifold defining function $g$ after an intermediate sub-step is given by
	\begin{equation}
		\tilde{g}_i \ppar*{v_{i+1}^{n[s+1]},t^{n[s+1]}} = \frac{\ppar*{1-\eta^{[s]}}^2}{2} \sum_{jk} \pdv{\tilde{g}_i}{v_j}{v_k} \ppar*{\zeta , t^{n[s+1]}} \cdot \ppar*{\Delta v_{i+1,j}^{n[s]} - \Delta v_{i,j}^{n[s]}} \cdot \ppar*{\Delta v_{i+1,k}^{n[s]} - \Delta v_{i,k}^{n[s]}}
	\end{equation}
	where $\zeta$ is some point in the orthotope with corners $v_i^{n[s+1]}$ and $v_{i+1}^{n[s+1]}$ and faces parallel to the axis.
\end{prop}
\begin{proof}
	By the Lagrange remainder theorem for Taylor series
	\begin{equation}\begin{split}
		g_i\ppar*{v_{i+1}^{n[s+1]},t^{n[s+1]}} ={}& g_i\ppar*{v_i^{n[s+1]},t^{n[s+1]}} + \sum_j \pdv{g_i}{v_j} \ppar*{v_i^{n[s+1]},t^{n[s+1]}} \cdot \ppar*{v_{i+1,j}^{n[s+1]}-v_{i,j}^{n[s+1]}} \\
		&+ \frac{1}{2} \sum_{ij} \pdv{g_i}{v_j}{v_k} \ppar*{\zeta,t^{n[s+1]}} \cdot \ppar*{v_{i+1,j}^{n[s+1]}-v_{i,j}^{n[s+1]}} \cdot \ppar*{v_{i+1,k}^{n[s+1]}-v_{i,k}^{n[s+1]}}.
	\end{split}\end{equation}
	Noting that
	\begin{equation}
		v_{i+1}^{n[s+1]}-v_i^{n[s+1]} = \ppar*{1-\eta^{[s]}} \cdot \ppar*{\Delta v_{i+1}^{n[s]} - \Delta v_{i}^{n[s]}}
	\end{equation}
	we next apply the scheme
	\begin{align}
		\begin{split}
		\pdv{\tilde{g}}{v} \ppar*{v_i^{n[s+1]},t^{n[s+1]}} \cdot \ppar*{1-\eta^{[s]}} \cdot \Delta v_{i+1}^{n[s]}
		={}& \ppar*{1-\eta^{[s]}} \cdot \tilde{G}\ppar*{v_i^{n[s+1]},t^{n[s+1]}} \cdot
		\pbrk*{ \begin{array}{c}
			\tilde{B}\ppar*{\Delta v_i^{n[s]} , t^{n[s]}}
			\\
			\ppar*{1-\eta^{[s]}} \cdot \tilde{G}\ppar*{v_i^{n[s+1]},t^{n[s+1]}}
		\end{array}}^{-1}
		\\
		& \cdot
		\pbrk*{ \begin{array}{c}
			\tilde{b}\ppar*{v^{n[s]},t^{n[s]}} \cdot \Delta t^n
			\\
			\ppar*{1-\eta^{[s]}} \cdot \tilde{G}\ppar*{v_i^{n[s+1]},t^{n[s+1]}} \cdot \Delta v_i^{n[s]} - \tilde{g}\ppar*{v_i^{n[s+1]},t^{n[s+1]}}
		\end{array}}
		\end{split}
		\\
		={}& \ppar*{1-\eta^{[s]}} \cdot \tilde{G}\ppar*{v_i^{n[s+1]},t^{n[s+1]}} \cdot \Delta v_i^{n[s]} - \tilde{g}\ppar*{v_i^{n[s+1]},t^{n[s+1]}},
	\end{align}
	and from these we obtain the result.
\end{proof}

Iterating the Newton-Raphson portion of the scheme yields, as expected,
\begin{equation}
	\norm*{\tilde{g}\ppar*{v_{i+1}^{n[s+1]},t^{n[s+1]}}} = \order{\norm*{\Delta v_{i+1}^{n[s]} - \Delta v_i^{n[s]}}^2}.
\end{equation}
This property of the scheme allows us to form an explicit version. If we initialise the Newton-Raphson iteration by setting $\Delta v_0^{n[s]}$ to the change in $v$ over the previous time-step then, in smooth regions, 
\begin{equation}
	\norm*{g\ppar*{v_1^{n[s+1]},t^{n[s+1]}}} = \order{\ppar*{\Delta t^n}^4}.
\end{equation}
This is more than sufficient considering there can be no drift-off error accumulation, and the overall Runge-Kutta scheme is, at most, third order. Following these observations we define the explicit version of this scheme to be
\begin{equation}
	\mathscr{E}_E\ppar*{v^{n[s]}} = \mathscr{E}\ppar*{v^{n[s]},\Delta v_0^{n[s]}}.
\end{equation}
We term these schemes RKNR (Runge-Kutta Newton-Raphson). In our test problems we will use the implicit version, but the presence of a highly accurate explicit modification indicates that the solution will converge rapidly.

%% file: Sections/Temporal_Comparison.tex
\subsection{Comparison of methods} \label{sec:TDP_tests}

\begin{figure}[tp!]
	\centering
	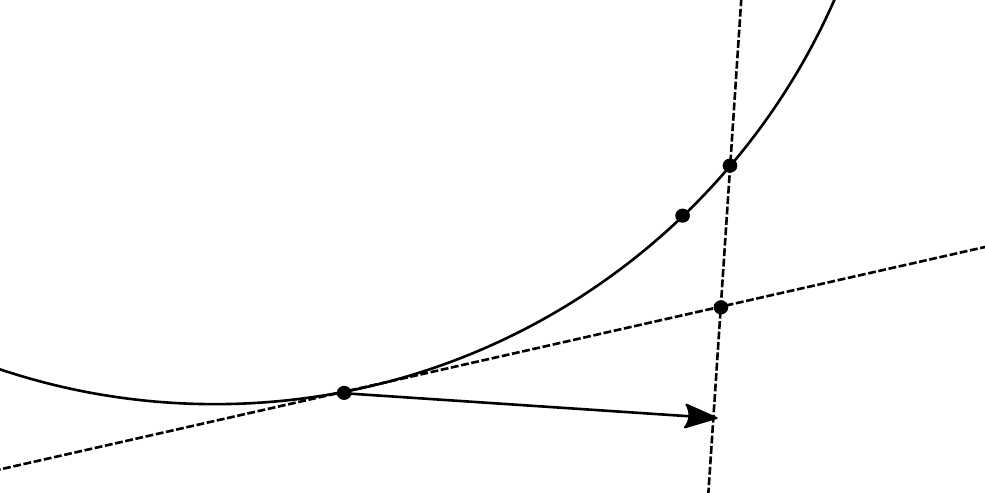
	\caption{Schematic representation of first order schemes in the case of dynamics in two-dimensional space confined to a one-dimensional time-invariant manifold with constant matrix $\tilde{B}$. The curve (a) represents the manifold, (b) the tangent at the starting point, and (c) the line defined by $\tilde{B} \Delta v = \tilde{b} \Delta t$, the value of $\Delta v$ parallel to $\tilde{B}^T$ is represented as an arrow. Marked are the points to which the different Runge-Kutta schemes send the solution, note that for first order schemes RKI1 is equivalent to RKP1 and RKI2 is equivalent to RKP2. If the local coordinate system on the manifold is an affine function of arc length then the distance between the starting point and the result of RK0 is equal to the arc length between the starting point and the result of RKLC. The angle formed at RKP1 between RK0 and the local tangent is only approximately a right angle (i.e. a right angle is obtained in the limit $\Delta t \rightarrow 0$), a true right angle would be obtained if using the optimisation stated in \cite{ln_Hairer_SDEM} directly.}
	\label{fig:TDP_schemes}
\end{figure}

The first order versions of the various schemes presented thus far are illustrated in \cref{fig:TDP_schemes}. We compare the basic Runge-Kutta scheme RK0 and the projection schemes RKP1, RKP2, RKI1, RKI2, and RKNR (all projected until the change in $v$ is less than $10^{-14}$) presented earlier in this section by testing them on a example problems. The class of problems considered is deliberately obtuse to distinguish between the different schemes, and are of the form
\begin{align}
	\frac{h_2(t)+10-y}{10}
	\ppar*{\begin{array}{cc}x&y\end{array}}
	\dv{}{t} \ppar*{\begin{array}{c}x\\y\end{array}}
	&= \dv{h_1}{t}
	&
	g(x,y,t) = \frac{h_2(t)+10-y}{10}(h_2(t) - y) = 0
\end{align}
were $h_1$ and $h_2$ are functions of $t$, and the factors of $h_2(t)+10-y$ increase the non-linearity of the system. The initial conditions are $x(t_0)=x_0$, $y(t_0)=y_0$ where $x_0^2 + y_0^2 = r_0^2$ and $y_0 = h_2(t_0)$, thus specifying $(r_0,t_0)$ identifies the solution
\begin{align}
	x(t) &= \sqrt{2\ppar*{h_1(t)-h_1(t_0)} + r_0^2 - h_2(t)^2},
	&
	y(t) &= h_2(t).
\end{align}
 In the test problems we will project the error into the row space of the differential equation using $\tilde{B}^T \ppar*{\tilde{B} \tilde{B}^T}^{-1} \tilde{B}$, the \emph{differential error}, and into the normal space of the manifold using $\tilde{G}^T \ppar*{\tilde{G} \tilde{G}^T}^{-1} \tilde{G}$, the \emph{algebraic error}. The error is computed from the projection by computing the absolute difference between the simulated and exact solutions after each time-step-projection, and then averaging over all time (\ie we use the $\ell_1$-norm).

Our test problems are all constructed on $0 \leq t \leq 1$ with $r_0 = 2$. For each test problem we run the simulation at a number of temporal resolutions, specifically
\begin{equation}	\label{eqn:TDP_resolutions}
	\setpred*{ 4 \round*{ \frac{10^{a/4}}{4} } }{ a \in \mathbb{N} \text{ and } 8 \leq a \leq 16 } = \pbrc*{100, 176, 316, 564, 1000, 1780, 3164, 5624, 10000}
\end{equation}
time-steps. Some of the test problems feature functions with discontinuous behaviour at $t \in \pbrc*{ 1/4 , 1/2 , 3/4 }$, which is why all of the resolutions have the same remainder on division by $4$.

\begin{testproblem}[analytic functions]\label{tp:TDP_analytic}

\begin{figure}[tp!]
	\centering
	\includegraphics{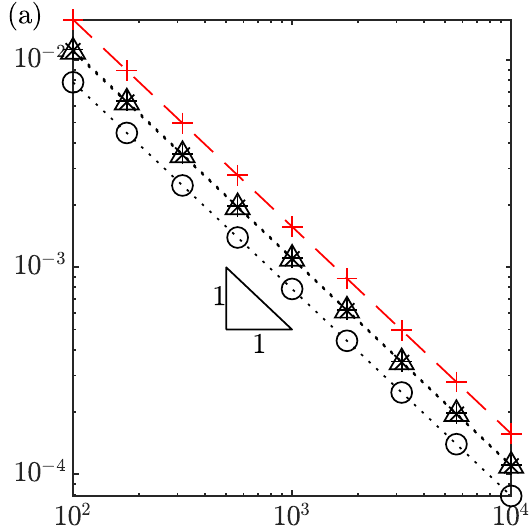}
	\includegraphics{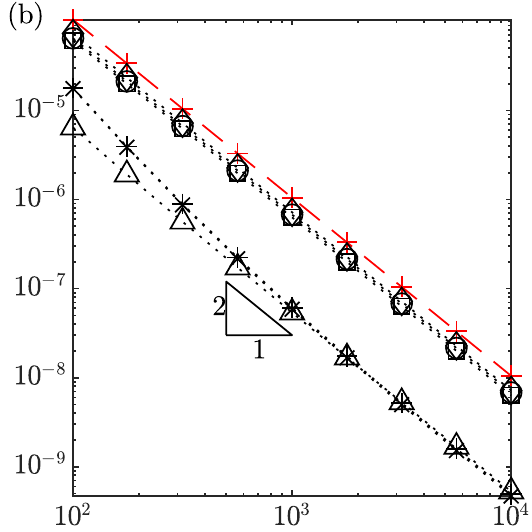}
	\includegraphics{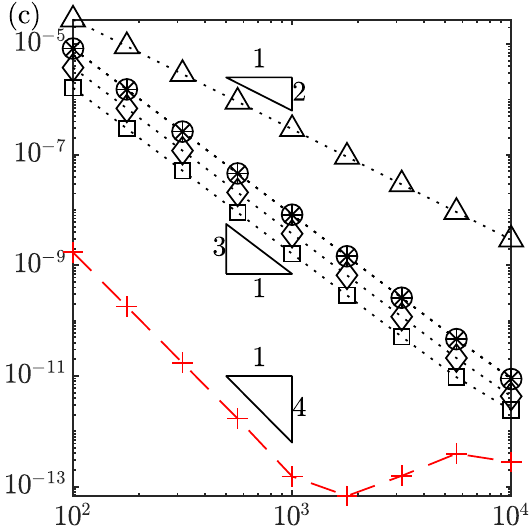}
	\caption{Convergence results for test problem \ref{tp:TDP_analytic}, in (a) we use $S=1$, (b) $S=2$, and (c) $S=3$, where $S$ is as in \cref{tab:RKcoeff}. In black with dotted lines is the differential error, while in red with dashed lines is the algebraic error (RK0 only), plotted as a function of the number of time-steps. The different projection schemes are distinguished by their markers: $+$ is RK0, $\times$ is RKP1, $\bigcirc$ is RKP2, $\square$ is RKI1, $\Diamond$ is RKI2, $\Delta$ is RKNR. We do not plot RKI1 or RKI2 for $S=1$, as these are the same as RKP1 and RKP2.}
	\label{fig:TDP_convergance_analytic}
\end{figure}

The first problem we consider is $h_1 = t^3/3$ and $h_2 = 1+\sin(2 \pi t)/2$. Thus the functions are all analytic, and we should obtain $S^\text{th}$ order convergence for a $S^\text{th}$ order scheme. Convergence results are plotted in \cref{fig:TDP_convergance_analytic}. We first note that the algebraic error for all projection schemes is zero, or around machine precision at $10^{-16}$. For RK0 the convergence of this error is at order $S$, as we would expect, except for $S=3$ which produces fourth order convergence (in fact, we have found that for $\tilde{g} = h_2(t)-y$, $h_2$ analytic, the use of (\ref{eqn:TDP_DAE_ODEversion}b) and the `third order' scheme produces fourth order convergence). For the differential error, at $S=1$ all schemes converge at first order, and are approximately of equal accuracy, with RKP2 being only slightly more accurate. At $S=2$ we see a more substantial divide, all schemes being second order, but RK0, RKP1, and RKNR are around an order of magnitude more accurate, RKNR having this edge even at low resolution. At $S=3$ we see some surprising behaviour, RKNR converging at second order, and at lower accuracy then the $S=2$ version. The other schemes are all third order, and have around the same accuracy, with RK and RKI2 being marginally better than the others.

\end{testproblem}

\begin{testproblem}[continuous functions]\label{tp:TDP_continuous}

\begin{figure}[tp!]
	\centering
	\includegraphics{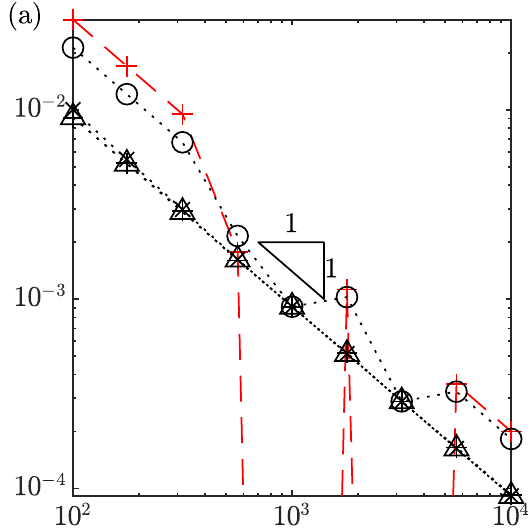}
	\includegraphics{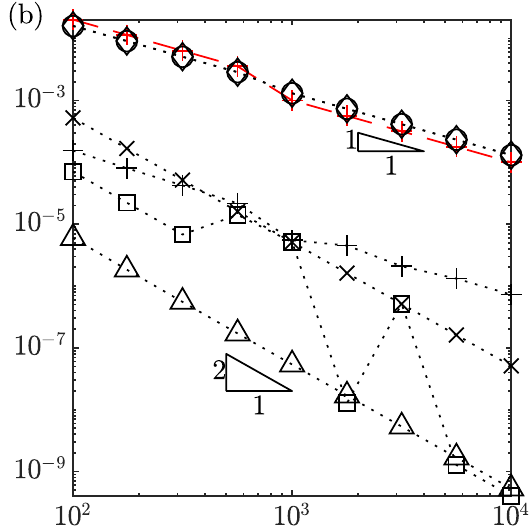}
	\includegraphics{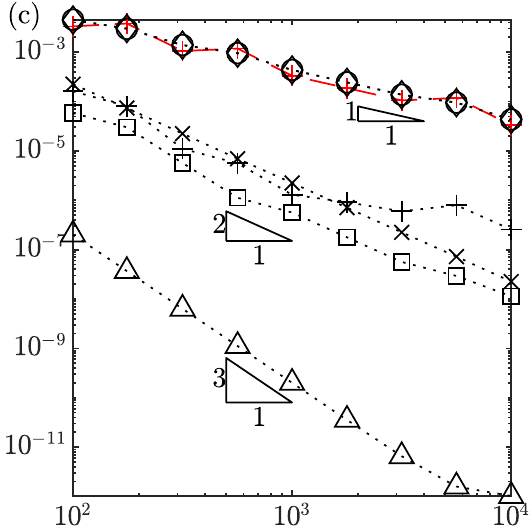}
	\caption{Convergence results for test problem \ref{tp:TDP_continuous}, in (a) we use $S=1$, (b) $S=2$, and (c) $S=3$, where $S$ is as in \cref{tab:RKcoeff}. In black with dotted lines is the differential error, while in red with dashed lines is the algebraic error (RK0 only), plotted as a function of the number of time-steps. The different projection schemes are distinguished by their markers: $+$ is RK0, $\times$ is RKP1, $\bigcirc$ is RKP2, $\square$ is RKI1, $\Diamond$ is RKI2, $\Delta$ is RKNR. We do not plot RKI1 or RKI2 for $S=1$, as these are the same as RKP1 and RKP2.}
	\label{fig:TDP_convergance_continuous}
\end{figure}

The next problem we consider is $h_1 = t^3/3$ and $h_2 = 1+\tri(t)/2$, where
\begin{align}
	\tri(t) &= 
	\begin{cases}
		4t & -1/4 \leq t \leq 1/4	\\
		2-4t & 1/4 < t < 3/4
	\end{cases},
	&
	\tri(t+1) &= \tri(t)
\end{align}
Thus $g(t)$, $v(t)$, and $\tilde{B}(v(t),t)$ have a discontinuous gradient. We expect all schemes to give at best first order convergence, the error being dominated by the effects near the discontinuous gradient. This is confirmed for the algebraic error for the RK0 scheme in \cref{fig:TDP_convergance_continuous}, with each increase in $S$ reducing the error by around half an order of magnitude. For the differential error, at $S=1$ we see first order convergence, with all schemes except RKP2 providing the same error, RKP2 varying depending on which side of the discontinuity expressions are evaluated. For $S=2$ there is a greater diversity, with RKNR being the most accurate, and RK0, RKP2 and RKI2 only converging at first order. At $S=3$, RKNR is the only scheme to achieve third order, which is even better than for \cref{tp:TDP_analytic}. We propose that this is because the algebraic condition satisfies $\pdv*[][3]{\tilde{g}}{t} = 0$ over most time-steps, which are therefore evaluated at third order, except local to the gradient discontinuities which introduce error over a single step.

\end{testproblem}

\begin{testproblem}[discontinuous functions]\label{tp:TDP_discontinuous}

\begin{figure}[tp!]
	\centering
	\includegraphics{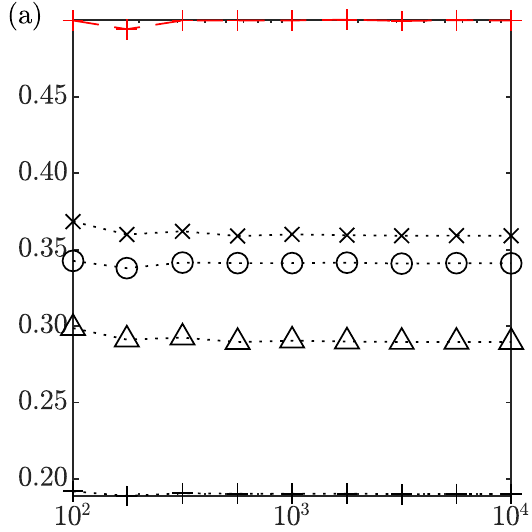}
	\includegraphics{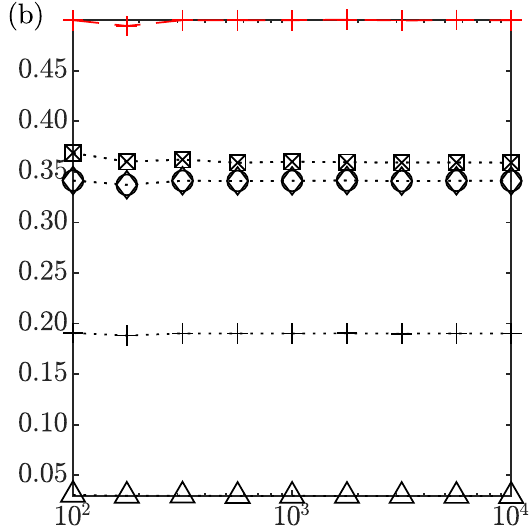}
	\includegraphics{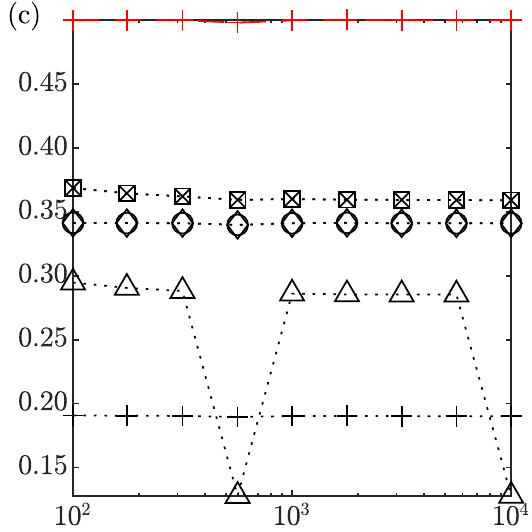}
	\caption{Convergence results for test problem \ref{tp:TDP_discontinuous}, in (a) we use $S=1$, (b) $S=2$, and (c) $S=3$, where $S$ is as in \cref{tab:RKcoeff}. In black with dotted lines is the differential error, while in red with dashed lines is the algebraic error (RK0 only), plotted as a function of the number of time-steps. The different projection schemes are distinguished by their markers: $+$ is RK0, $\times$ is RKP1, $\bigcirc$ is RKP2, $\square$ is RKI1, $\Diamond$ is RKI2, $\Delta$ is RKNR. We do not plot RKI1 or RKI2 for $S=1$, as these are the same as RKP1 and RKP2.}
	\label{fig:TDP_convergance_discontinuous}
\end{figure}

We consider $h_1 = t^3/3$ and $h_2 = 1+\sqr(t)/2$, where
\begin{align}
	\sqr(t) &= 
	\begin{cases}
		1 & 0 \leq t < 1/2	\\
		-1 & 1/2 \leq t < 1
	\end{cases},
	&
	\sqr(t+1) &= \sqr(t).
\end{align}
Thus $g(t)$, $v(t)$, and $\tilde{B}(v(t),t)$ have a discontinuity. We should expect none of the schemes to converge, the error being dominated by the effects near the discontinuity. \Cref{fig:TDP_convergance_discontinuous} shows us exactly this, the algebraic error in RK0 being largely independent of $S$ and resolution. The differential error at $S=1$ shows that RK0 is the most accurate, though this is largely irrelevant as it cannot resolve the manifold. Of the other schemes RKNR has the least error at all orders.

\end{testproblem}

\begin{testproblem}[singular forcing]\label{tp:TDP_singular}

\begin{figure}[tp!]
	\centering
	\includegraphics{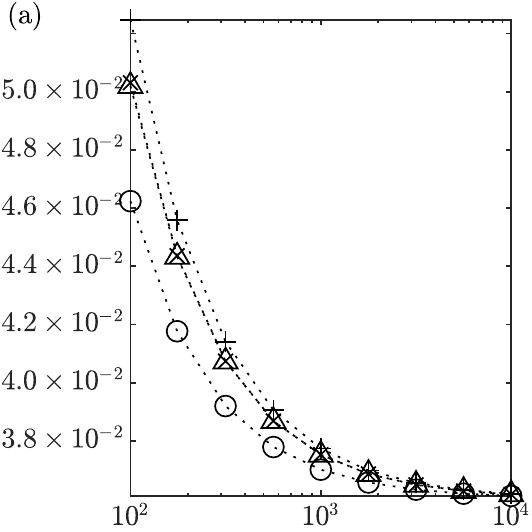}
	\includegraphics{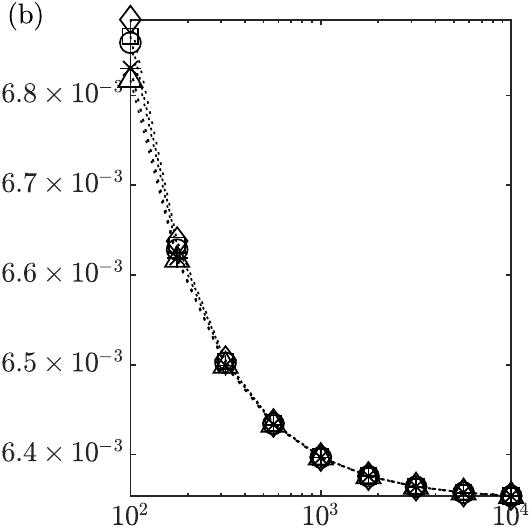}
	\includegraphics{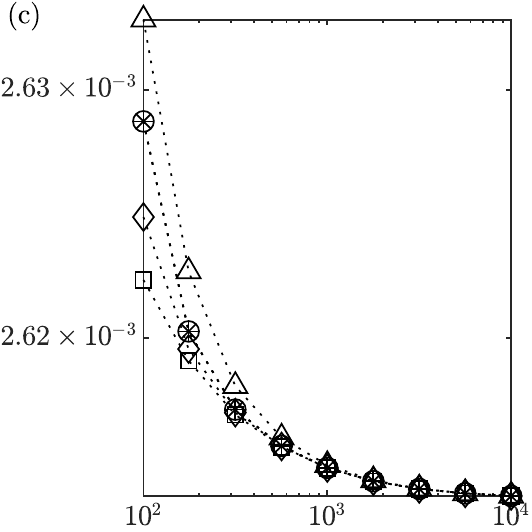}
	\caption{Convergence results for test problem \ref{tp:TDP_singular}, in (a) we use $S=1$, (b) $S=2$, and (c) $S=3$, where $S$ is as in \cref{tab:RKcoeff}. In black with dotted lines is the differential error plotted as a function of the number of time-steps. The different projection schemes are distinguished by their markers: $+$ is RK0, $\times$ is RKP1, $\bigcirc$ is RKP2, $\square$ is RKI1, $\Diamond$ is RKI2, $\Delta$ is RKNR. We do not plot RKI1 or RKI2 for $S=1$, as these are the same as RKP1 and RKP2. The algebraic error is not plotted, and is approximately the same as in \cref{fig:TDP_convergance_analytic}.}
	\label{fig:TDP_convergance_singular}
\end{figure}

The final problem we consider is singular forcing of the ODE over a time period proportional to $\Delta t$. We take $h_1 = \tanh((t-1/2)/\Delta t)$ and $h_2 = 1+\sin(2 \pi t)/2$. We expect that error generation will be primarily around $t=1/2$ where the forcing occurs, and as this region scales with $\Delta t$ that the error should be invariant of the resolution. Plotted in \cref{fig:TDP_convergance_singular} is the differential error, for which we do indeed see this behaviour. The different projection methods have little effect because the solution never moves far from the manifold. The algebraic error for RK0 is found to be around the same as for \cref{tp:TDP_analytic}, which is consistent with both cases having the same algebraic equation. It seems then that singular forcing has little impact on the convergence to the manifold.

\end{testproblem}

\begin{discuss}

The results of these test problems tell us that, in smooth regions, first order schemes are all largely similar, of the second order schemes we should use RKP1 or RKNR, and at third order we should use the RKI1 or RKI2, as these are the most accurate of the schemes for test problem \ref{tp:TDP_analytic}. However, not all problems for hyperbolic systems have smooth solutions, many contain shocks which are diffused over a small number of grid cells to produce a rapid transition as investigated in test problem \ref{tp:TDP_singular}. There we see that the choice of projection scheme is not particularly relevant as all schemes generate approximately the same error, which does not impact the convergence to the manifold. Some problems, though, include sudden changes in algebraic boundary conditions, the boundary conditions being defined piecewise to enforce different conditions depending on the nature of the flow in the bulk. It is very desirable to minimise the amount of error generated during a change in condition, as illustrated in \cref{tp:TDP_continuous,tp:TDP_discontinuous}. If a second order scheme is sufficient then RKNR is much better than the other schemes for both discontinuous gradients and slightly better for discontinuities. If a third order scheme is required then RKP1 or RKI1 are better than the others for discontinuous gradients, and similar to the rest for discontinuities, and considering computational efficiency we therefore recommend RKP1. The $S=3$ RKNR is not recommended because it cannot be guaranteed to be third order by \cref{tp:TDP_analytic}. All schemes are quite poor around discontinuities in algebraic boundary condition, see \cref{tp:TDP_discontinuous}, and as a result event detection may be required so that a carefully designed treatment to problems that arise may be applied before continuing simulation.

We observe that the RKI1, RKI2, and RKNR scheme can be more computationally intensive than the post step projection schemes RKP1 and RKP2. However, for simulating \cref{eqn:hyp_system} when the bulk is evaluated at high resolution this dominates the computational time so implicit evolution of the boundary conditions is an excusable expense.

As a final note, even though in these test problems we obtained $\order{1}$ errors, we do not expect this to carry over to the simulation of hyperbolic systems. This is because, for outgoing characteristics, the solution will be restored to the values contained in the bulk scheme over a timescale $\order{\Delta t}$, so we expect first order convergence. Errors in the algebraic conditions, static characteristics, and non-reflecting conditions are not restored in this manner as their values are advected into the domain, which is why we require them to be evolved at higher accuracy.

\end{discuss}

%% file: Figures/RKNRstepping.pdf_tex
\begingroup%
  \makeatletter%
  \providecommand\color[2][]{%
    \errmessage{(Inkscape) Color is used for the text in Inkscape, but the package 'color.sty' is not loaded}%
    \renewcommand\color[2][]{}%
  }%
  \providecommand\transparent[1]{%
    \errmessage{(Inkscape) Transparency is used (non-zero) for the text in Inkscape, but the package 'transparent.sty' is not loaded}%
    \renewcommand\transparent[1]{}%
  }%
  \providecommand\rotatebox[2]{#2}%
  \newcommand*\fsize{\dimexpr\f@size pt\relax}%
  \newcommand*\lineheight[1]{\fontsize{\fsize}{#1\fsize}\selectfont}%
  \ifx\svgwidth\undefined%
    \setlength{\unitlength}{283.46456693bp}%
    \ifx\svgscale\undefined%
      \relax%
    \else%
      \setlength{\unitlength}{\unitlength * \real{\svgscale}}%
    \fi%
  \else%
    \setlength{\unitlength}{\svgwidth}%
  \fi%
  \global\let\svgwidth\undefined%
  \global\let\svgscale\undefined%
  \makeatother%
  \begin{picture}(1,0.5)%
    \lineheight{1}%
    \setlength\tabcolsep{0pt}%
    \put(0,0){\includegraphics[width=\unitlength,page=1]{RKNRstepping.pdf}}%
    \put(0.84400574,0.45944351){\color[rgb]{0,0,0}\makebox(0,0)[lt]{\lineheight{1.25}\smash{\begin{tabular}[t]{l}(a)\end{tabular}}}}%
    \put(0.9394149,0.26196106){\color[rgb]{0,0,0}\makebox(0,0)[lt]{\lineheight{1.25}\smash{\begin{tabular}[t]{l}(b)\end{tabular}}}}%
    \put(0.72560661,0.00944366){\color[rgb]{0,0,0}\makebox(0,0)[lt]{\lineheight{1.25}\smash{\begin{tabular}[t]{l}(c)\end{tabular}}}}%
    \put(0.74579468,0.35638428){\color[rgb]{0,0,0}\makebox(0,0)[rt]{\lineheight{1.25}\smash{\begin{tabular}[t]{r}RKNR \& RKP2\end{tabular}}}}%
    \put(0.69524803,0.29670197){\color[rgb]{0,0,0}\makebox(0,0)[rt]{\lineheight{1.25}\smash{\begin{tabular}[t]{r}RKLC\end{tabular}}}}%
    \put(0.73874489,0.14917572){\color[rgb]{0,0,0}\makebox(0,0)[lt]{\lineheight{1.25}\smash{\begin{tabular}[t]{l}RK0\end{tabular}}}}%
    \put(0,0){\includegraphics[width=\unitlength,page=2]{RKNRstepping.pdf}}%
    \put(0.64514145,0.24566467){\color[rgb]{0,0,0}\makebox(0,0)[rt]{\lineheight{1.25}\smash{\begin{tabular}[t]{r}RKP1\end{tabular}}}}%
  \end{picture}%
\endgroup%

%% file: Sections/Test.tex
\section{Implementation and numerical tests for the shallow water equations} \label{sec:FVTests}

In this section we will present results from simulations of the shallow water equations, which we use as a test case for our numerical implementation of boundary conditions. Prior to this we detail the numerical scheme.

\input{Sections/Test_RescaleDomain}
\input{Sections/Test_ShallowWater}
\input{Sections/Test_FiniteVolume}

\input{Sections/Test_Problems}

%% file: Sections/Test_RescaleDomain.tex
\subsection{Temporally evolving domain}	\label{sec:TED}

We begin by transforming from the temporally evolving domain $x \in [x_L(t),x_R(t)]$ to a fixed domain $\mathscr{y} \in [0,1]$ using the transformation (\eg \cite{bk_Ungarish_GCI})
\begin{subequations}	\label{eqn:TED_transform}
	\begin{align}
	\mathscr{y} &= \mathscr{Y}(x,t) \eqdef \frac{x-x_L(t)}{x_R(t)-x_L(t)}
	& \text{and}&&
	x &= \mathscr{X}(\mathscr{y},t) \eqdef (1-\mathscr{y}) \cdot x_L(t) + \mathscr{y} \cdot x_R(t).
\end{align}
We define special notation for the derivatives of $\mathscr{X}$, namely
\begin{align}	\label{eqn:TED_rl_def}
	\mathscr{r}(\mathscr{y},t)\eqdef \pdv{\mathscr{X}}{t} &= (1-\mathscr{y}) \cdot \dv{x_L}{t} (t) + \mathscr{y} \cdot \dv{x_R}{t} (t),
	&\text{and}&&
	\mathscr{l}(t) \eqdef \pdv{\mathscr{X}}{\mathscr{y}} &= x_R(t) - x_L(t),
	&
	\pgrp*{ \pdv{\mathscr{r}}{\mathscr{y}} \equiv \dv{\mathscr{l}}{t} }.
\end{align}
\end{subequations}
The variable $\mathscr{r}$ we call the \emph{rate} and is the speed of the curves of constant $\mathscr{y}$ through $(x,t)$ space, the variable $\mathscr{l}$ is the domain \emph{length}. Using this transformation, ODEs transform from
\begin{align}\label{eqn:transform_speed}
	\dv{x}{t} &= \lambda(x,t)
	&\text{to}&&
	\dv{\mathscr{y}}{t} &= \frac{\lambda(\mathscr{X}(\mathscr{y},t),t) - \mathscr{r}(\mathscr{y},t)}{\mathscr{l}(t)},
\end{align}
which tells us how the characteristic speeds $\lambda^{(m)}$ transform.

To transform the balance law \cref{eqn:hyp_system} we first write it in the form 
\begin{equation}
	\pdv{Q}{t} + \pdv{}{x} \ppar*{\mathscr{u}Q + \hat{F}} = \Psi
\end{equation}
where $\mathscr{u}:(Q,x,t) \mapsto \mathbb{R}$ is a velocity representing how fast we expect $Q$ to be advected; if a physical velocity is present we may use this for $\mathscr{u}$. We could simply transform our system using \eqref{eqn:TED_transform}, but this has been found to introduce additional source terms, which is undesirable for systems where $\Psi \equiv 0$. Instead we introduce transformations of the density, additional flux $\hat{F} \eqdef F - \mathscr{u} Q$, and source as
\begin{subequations}	\label{eqn:TED_transform_varables}
\begin{align}
	\mathscr{Q}(\mathscr{y},t) &= \mathscr{T}\ppar*{\mathscr{r}(\mathscr{y},t),\mathscr{l}(t)} \cdot Q\ppar*{\mathscr{X}(\mathscr{y},t) , \: t }, \label{eqn:TED_transform_varables_Q} 
	\\
	\hat{\mathscr{F}}(\mathscr{Q},\mathscr{y},t) &= \mathscr{T}\ppar*{\mathscr{r}(\mathscr{y},t),\mathscr{l}(t)} \cdot \hat{F}\ppar*{\mathscr{T}\ppar*{\mathscr{r}(\mathscr{y},t),\mathscr{l}(t)}^{-1}\mathscr{Q} , \: \mathscr{X}(\mathscr{y},t) , \: t }, 
	\\
	\mathscr{\hat{S}}(\mathscr{Q},\mathscr{y},t) &= \mathscr{T}\ppar*{\mathscr{r}(\mathscr{y},t),\mathscr{l}(t)} \cdot \Psi\ppar*{\mathscr{T}\ppar*{\mathscr{r}(\mathscr{y},t),\mathscr{l}(t)}^{-1}\mathscr{Q} , \: \mathscr{X}(\mathscr{y},t) , \: t },
\end{align}
\end{subequations}
where $\mathscr{T}:(\mathscr{r},\mathscr{l}) \mapsto \mathbb{M}(M,M)$. We require that $\mathscr{\mathscr{T}}$ be continuous as a function and invertible as a matrix so that the weak solutions to the system transform correctly. Making this change of variables yields the system
\begin{align} \label{eqn:hyp_system_transform}
	\pdv{\mathscr{Q}}{t} +\pdv{}{\mathscr{y}}\ppar*{ \mathscr{F} } = \mathscr{S}
\end{align}
where
\begin{align}	\label{eqn:TED_resultant_flux_source}
	\mathscr{F} &\eqdef \frac{(\mathscr{u}-\mathscr{r})\mathscr{Q} + \hat{\mathscr{F}}}{\mathscr{l}},
	&
	\mathscr{S} &\eqdef  \pdv{\mathscr{r}}{t}\pdv{\mathscr{T}}{\mathscr{r}} \mathscr{T}^{-1} \mathscr{Q} + \dv{\mathscr{l}}{t} \ppar*{ \pdv{\mathscr{T}}{\mathscr{l}} - \frac{\mathscr{T}}{\mathscr{l}} + \frac{\mathscr{u}-\mathscr{r}}{\mathscr{l}} \pdv{\mathscr{T}}{\mathscr{r}} } \mathscr{T}^{-1} \mathscr{Q} + \frac{1}{\mathscr{l}}\dv{\mathscr{l}}{t}\pdv{\mathscr{T}}{\mathscr{r}} \mathscr{T}^{-1} \hat{\mathscr{F}} + \hat{\mathscr{S}}.
\end{align}
While this equation looks far worse than the original, we have freedom to choose $\mathscr{T}$ to provide properties that we desire. For example, if we simply wish to eliminate as many source terms as possible then we can take $\mathscr{T} = \mathscr{l} I$ where $I \in \mathbb{M}(M,M)$ is the identity, which yields the source $\mathscr{S} =  \hat{\mathscr{S}}$. Alternatively we may be interested in ensuring the $m^\text{th}$ field is conserved to machine precision, which may be achieved by constructing $\mathscr{T}$ such that $\mathscr{F}_{m}=0$ at the boundaries and $\mathscr{S}_{m} = 0$ everywhere. We discuss how this may be done for the shallow water equations next.

%% file: Sections/Test_ShallowWater.tex
\subsection{The shallow water system} \label{sec:SW}

We demonstrate our numerical scheme by application to a shallow-water model of fluid flow containing a passive tracer. This models an open channel flow driven by the difference in fluid densities between the current and surrounding dynamically passive ambient. We denote by $h(x,t)$ the fluid depth, $u(x,t)$ the fluid velocity, and $\phi(x,t)$ the volumetric concentration of tracer. The (dimensionless) system is (\eg \cite{bk_Ungarish_GCI})
\begin{align}
	\pdv{h}{t} + \pdv{}{x}\ppar*{uh} &= 0,
	&
	\pdv{}{t}(\phi h) + \pdv{}{x}\ppar*{u \phi h} &= 0,
	&
	\pdv{}{t}(uh) + \pdv{}{x}\ppar*{u^2 h + \frac{h^2}{2}} &= 0.
\end{align}
We transform this system from the domain $x_L(t) \leq x \leq x_R(t)$ to $0 \leq \mathscr{y} \leq 1$ as discussed in \cref{sec:TED}. An important boundary condition is $u=\mathscr{r}$, as this enforces no fluid flux through the domain end, and the ability to simulate this exactly means that, in a simulation with this boundary condition at both ends of the domain, the volume of fluid will be conserved to machine precision. To achieve this we take $\mathscr{u} \eqdef u$, and our transformed variables are
\begin{align}
	\whatchar{h} &\eqdef f_1(\mathscr{l}) \cdot h, &
	\widehat{\phi h} &\eqdef f_2(\mathscr{l}) \cdot \phi h, &
	\widehat{u h} &\eqdef f_3(\mathscr{l}) \cdot (u-\mathscr{r}) \cdot h,
\end{align}
for functions $f_i$ to be determined. The transformed density is then
\begin{align}
	\mathscr{Q} &\eqdef
	\pgrp*{\begin{array}{c}
		\whatchar{h}			\\
		\widehat{\phi h}	\\
		\widehat{u h}
	\end{array}}
	=
	\mathscr{T} Q
	&\text{where}&&
	\mathscr{T} &\eqdef
	\pgrp*{\begin{array}{ccc}
		f_1(\mathscr{l}) 				&	0						&	0		\\
		0 								&	f_2(\mathscr{l})	&	0		\\
		-\mathscr{r} f_3(\mathscr{l}) 	&	0						&	f_3(\mathscr{l})
	\end{array}},
	&
	Q &\eqdef
	\pgrp*{\begin{array}{c}
		h	\\
		\phi h	\\
		u h
	\end{array}}.
\end{align}
We see from \eqref{eqn:TED_resultant_flux_source} that we can remove a large portion of the resultant source terms by requiring that
\begin{align}
	0
	= \ppar*{\pdv{\mathscr{T}}{\mathscr{l}} - \frac{\mathscr{T}}{\mathscr{l}} + \frac{u-\mathscr{r}}{\mathscr{l}} \pdv{\mathscr{T}}{\mathscr{r}}} Q
	&=
	\pgrp*{\begin{array}{ccc}
		f_1' - f_1/\mathscr{l}					&	0						&	0					\\
		0 										&	f_2' - f_2/\mathscr{l}	&	0					\\
		-\mathscr{r} (f_3' - 2 f_3/\mathscr{l})	&	0						&	f_3' - 2 f_3/\mathscr{l}
	\end{array}}
	\pgrp*{\begin{array}{c}
		h	\\
		\phi h	\\
		u h
	\end{array}},
\end{align}
where prime denotes $\dv*{}{\mathscr{l}}$. We take $f_1(\mathscr{l}) = f_2(\mathscr{l})= \mathscr{l}$, and $f_3(\mathscr{l}) = \mathscr{l}^2$, and the resultant system is
\begin{align}
	\pdv{}{t} \ppar*{\whatchar{h}} + \pdv{}{y} \ppar*{ \frac{\widehat{uh}}{\mathscr{l}^2} } &= 0,
	&
	\pdv{}{t} \ppar*{\widehat{\phi h}} + \pdv{}{y} \ppar*{ \frac{\widehat{uh}\widehat{\phi h}}{\mathscr{l}^2 \whatchar{h}} } &= 0,
	&
	\pdv{}{t} \ppar*{\widehat{u h}} + \pdv{}{y} \ppar*{ \frac{\widehat{u h}^2}{\mathscr{l}^2 \whatchar{h}} + \frac{\whatchar{h}^2}{2\mathscr{l}} } &= - \mathscr{l} \whatchar{h}\pdv{\mathscr{r}}{t}.
\end{align}
The characteristic speeds and left eigenvectors of this system are
\begin{subequations} \label{eqn:SW_eigen}
\begin{align} 
	\label{eqn:SW_eigenvalue}
	\whatchar{\lambda}^{(1)}	&= \frac{\widehat{uh}}{\mathscr{l}^2\whatchar{h}} - \frac{\sqrt{\whatchar{h}}}{\mathscr{l}^{3/2}}, &
	\whatchar{\lambda}^{(2)}	&= \frac{\widehat{uh}}{\mathscr{l}^2\whatchar{h}},	&
	\whatchar{\lambda}^{(3)}	&= \frac{\widehat{uh}}{\mathscr{l}^2\whatchar{h}} + \frac{\sqrt{\whatchar{h}}}{\mathscr{l}^{3/2}},
\\
	\whatchar{l}^{(1)}	&= 
	\pgrp*{\begin{array}{ccc}
		\mathscr{l}^2 \whatchar{\lambda}^{(3)},	&
		0, 	&
		- 1
	\end{array}},
	&
	\whatchar{l}^{(2)}	&= 
	\pgrp*{\begin{array}{ccc}
		\widehat{\phi h},	&
		- \whatchar{h}, 	&
		0
	\end{array}},
	&
	\whatchar{l}^{(3)}	&= 
	\pgrp*{\begin{array}{ccc}
		\mathscr{l}^2 \whatchar{\lambda}^{(1)},	&
		0, 	&
		- 1
	\end{array}}.
\end{align}
\end{subequations}
In regions with continuous solution and non-zero depth, this system of equations has invariants on its characteristics
\begin{subequations}\begin{align}
	\alpha &\eqdef u + 2\sqrt{h} 
	&\text{is constant on any curve $\mathscr{y}(t)$ satisfying}&&
	\dv*{\mathscr{y}}{t} = \whatchar{\lambda}^{(3)},
	\\
	\beta &\eqdef u - 2\sqrt{h} 
	&\text{is constant on any curve $\mathscr{y}(t)$ satisfying}&&
	\dv*{\mathscr{y}}{t} = \whatchar{\lambda}^{(1)},
	\\
	\phi & 
	&\text{is constant on any curve $\mathscr{y}(t)$ satisfying}&&
	\dv*{\mathscr{y}}{t} = \whatchar{\lambda}^{(2)}.
\end{align}\end{subequations}

%% file: Sections/Test_FiniteVolume.tex
\subsection{Finite volume discretization} \label{sec:FV}

Once we have performed the desired transformations, the hyperbolic system \cref{eqn:hyp_system} will be in the general form \cref{eqn:hyp_system_transform} on $0 \leq \mathscr{y} \leq 1$. These equations have exactly the same structure, thus all of our discussion of boundary conditions is applicable to \cref{eqn:hyp_system_transform}. The system is spatially discretized over $J$ cells by introducing $J+1$ cell interfaces $0=\mathscr{y}_{1/2}<\mathscr{y}_{3/2}<\ldots<\mathscr{y}_{J+1/2}=1$, and defining cell averages
\begin{align}
	\mathscr{Q}_{j}(t) &= \frac{1}{\Delta \mathscr{y}_j} \int_{\mathscr{y}_{j-1/2}}^{\mathscr{y}_{j+1/2}} \mathscr{Q} \dd y,
	&
	\mathscr{S}_{j}(t) &= \frac{1}{\Delta \mathscr{y}_j} \int_{\mathscr{y}_{j-1/2}}^{\mathscr{y}_{j+1/2}} \mathscr{S} \dd y,
\end{align}
where $\Delta \mathscr{y}_j = \mathscr{y}_{j+1/2} - \mathscr{y}_{j-1/2}$ is the width of the $j^\text{th}$ cell. Averaging the system \eqref{eqn:hyp_system_transform} yields
\begin{equation}
	\dv{\mathscr{Q}_{j}}{t} + \frac{1}{\Delta \mathscr{y}_j} \pbrk*{ \mathscr{F}_{j+1/2} - \mathscr{F}_{j-1/2} } = \mathscr{S}_{j}
\end{equation}
where $\mathscr{F}_{j+1/2}(t) = \mathscr{F}(\mathscr{Q}(\mathscr{y}_{j+1/2},t),\mathscr{y}_{j+1/2},t)$. To construct a spatially discrete scheme we approximate $\mathscr{F}_{j+1/2}$ and $\mathscr{S}_{j}$ in terms of the discrete solution, which can be done using a variety of methods, \eg \cite{bk_Leveque_FVM,ar_Pirozzoli_2011}. We employ the the central upwind scheme \cite{ar_Kurganov_2001} for which $C_{\max} = 1/8$, and second order RKNR time-stepping unless otherwise stated.

To simulate the boundary conditions we use the method developed in \cref{sec:spatial} employing first order extrapolation. For demonstration purposes we calculate the flux directly from the computed values at the boundaries, \ie $\mathscr{F}_{J+1/2} = \mathscr{F}(\mathscr{Q}_R,1,t)$ (in applications it may be desirable impose $\mathscr{Q}_R$ as a Dirichlet condition using a specialist technique as discussed in \cref{sec:Intro}). The forcing will usually only be applied to outgoing and static characteristics, that is
\begin{subequations}
\begin{align} \label{eqn:diffusion_outstatic}
 	D^{(m)} = \frac{1}{\Delta \mathscr{y}_{J+1/2}}
 	\begin{cases}
 		4 \mathcal{D} \whatchar{\lambda}_{R\max} 	& \text{for } \max \ppar*{\whatchar{\lambda}_R^{(m)} , \whatchar{\lambda}_J^{(m)}} > -\delta,	\\
 		0 													& \text{otherwise}.
 	\end{cases}
\end{align}
where for our implementable we take $\delta = 10^{-8}$, but for some simulations in \cref{tp:TP_alphawave} we will apply forcing to incoming characteristics also, in which case
\begin{align} \label{eqn:diffusion_all}
	D^{(m)} = \frac{4 \mathcal{D} \whatchar{\lambda}_{R\max}}{ \Delta \mathscr{y}_{J+1/2}},
\end{align}
where we use (in \cref{tp:TP_critwall} this is modified to aid discussion of eigenvalues and eigenvectors)
\begin{align}
	\whatchar{\lambda}_{R\max} = \max \ppar*{ \max_{m} \abs*{\whatchar{\lambda}^{(m)}_R} , \max_{m} \abs*{\whatchar{\lambda}^{(m)}_J} }
\end{align}
\end{subequations}
and $\lambda^{(m)}_R$ are the values at $\mathscr{y} = 1$ while $\lambda^{(m)}_J$ are those at $\mathscr{y} = \mathscr{y}_J$. In the assembling of the discretized \cref{eqn:BC_extp_system_general,eqn:BC_extp_system_reflect} we operate in the transformed variables and use the eigenvalues and eigenvectors computed at $\mathscr{y} = 1$ as discussed in \cref{tp:TP_critwall}. Further details of the implementation may be found in \cref{sec:IMP}.

%% file: Sections/Test_Problems.tex
\subsection{Test problems}	\label{sec:TP}

In this section we document simulations of the shallow water system in order to demonstrate the effectiveness of our method. We will use the $\ell_1$ error for our convergence analysis, which is computed by taking the absolute difference between simulated and exact values for every variable and location stated, and then averaging over all.

We begin with \cref{tp:TP_critwall,tp:TP_particle} which expose the non-hyperbolic degeneracy of the system of equations when $h \rightarrow 0$, as well issues that occur when switching between boundary conditions. The computation of eigenvalues and eigenvectors is discussed, and the necessity of extrapolation and forcing (especially for static characteristics). In \cref{tp:TP_alphawave} we validate our stability analysis, and discuss non-reflecting conditions. In \cref{tp:TP_lockrel} we demonstrate the $\order{1}$ error generation that can occur when using RK0 time-stepping, and that projection is capable of remedying it. For all test problems we discuss convergence, using uniform spatial discretization with resolutions
\begin{align}
	J \in \setpred*{ \round*{ 10^{a/4} } }{ a \in \mathbb{N} \text{ and } 8 \leq a \leq 16 } = \pbrc*{ 100, 178, 316, 562, 1000, 1778, 3162, 5623, 10000 }.
\end{align}

\begin{testproblem}[overtopping a barrier] \label{tp:TP_critwall}

We begin with the consideration of a dam-break colliding with a barrier of finite height. The domain is $0 \leq x \leq 1$ (thus $\mathscr{y} = x$ and $\mathscr{Q} = Q$), and the initial condition is $h=\phi=1$ and $u=0$ on $0 \leq x \leq 1/2$, and $h=0$ on $1/2 < x \leq 1$. At  $x=0$ we enforce $u=0$ to prevent outflow, while at $x=1$ the dynamics are much more involved. As discussed in \cite{ar_Skevington_F001_Draining} the energy difference between the flow at the bottom of the barrier ($x=1$) and the minimal energy of flow at the top ($x=1^+$, immediately beyond the end of the domain), $\Delta E$, is (\cf \cref{eqn:BC_example})
\begin{align}
	\Delta E &= \eval*{ \frac{u^2}{2} + (h-h_b) - \frac{3}{2} ( uh )^{2/3} }_{\mathrlap{x = 1}},
\end{align}
where $h_b$ is the dimensionless height of the barrier; we take $h_b=1/2$. If the flow at the bottom is supercritical, $\lambda^{(1)}|_{x=1} > 0$, and its energy is sufficient, $\Delta E > 0$, then no boundary condition is required at $x=1$. Otherwise we enforce $\Delta E = 0$ so long as this is possible given the information flowing out from the bulk (\ie there is a solution to the corresponding $\Delta E = 0$, $u|_{x=1} > 0$, $\lambda^{(1)}|_{x=1} < 0$ boundary Riemann problem). If this is not possible then enforce no flow, $u=0$; the flow is too shallow/slow to overtop the barrier. Numerical details are given in \cref{sec:IMP}.

\begin{figure}[tp!]
	\centering
	\includegraphics{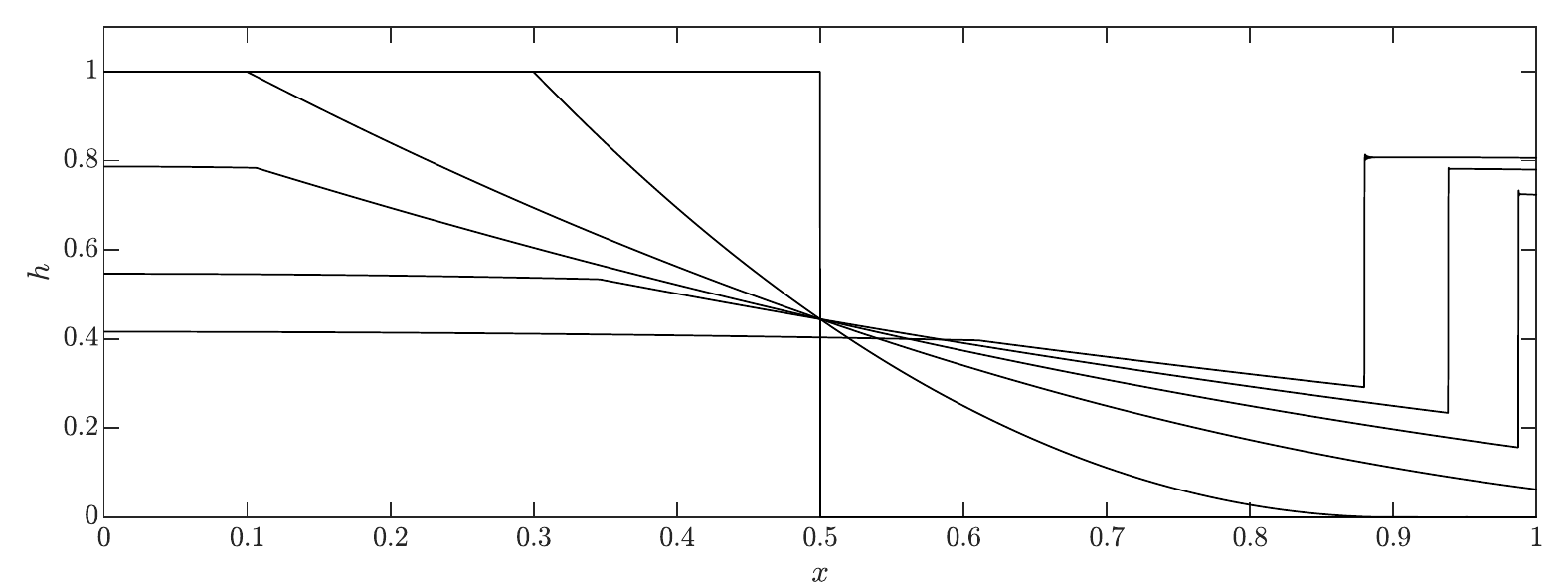}
	\caption{The depth at times $t\in \pbrc*{0, 0.2, \ldots , 1}$ for \cref{tp:TP_critwall} simulated at resolution $J=10^4$ with $\mathcal{D}=1/2$.}
	\label{fig:TP_critwall}
\end{figure}

The behaviour given these conditions, as shown in \cref{fig:TP_critwall}, is that the fluid to undergoes dam-break and wet the region between $x=1/2$ and $1$ as a rarefaction fan, \ie
\begin{align}\label{eqn:TP_critwall_exact}
	h &= 
	\begin{cases}
		1 								& \text{for }  x-1/2 \leq -t,			\\
		\tfrac{1}{9}(2-(x-1/2)/t)^2		& \text{for } -t \leq x-1/2 \leq 2t,	\\
		0 								& \text{for } 2t \leq x-1/2,
	\end{cases}
	&
	u &= 
	\begin{cases}
		0 								& \text{for }  x-1/2 \leq -t,			\\
		\tfrac{2}{3}(1+(x-1/2)/t)		& \text{for } -t \leq x-1/2 \leq 2t,	\\
		\text{undefined} 				& \text{for } 2t \leq x-1/2.
	\end{cases}
\end{align}
When the fluid arrives at the barrier ($t=1/4$) it is both supercritical and highly energetic, easily undergoing supercritical overtopping, thus no boundary condition is required and \cref{eqn:TP_critwall_exact} is still valid. As time passes, $\Delta E$ reduces until eventually we must enforce no flow at time $t=0.54$ (2 \sf), generating a backwards propagating shock. The fluid between the shock and the barrier rapidly deepens and subcritical outflow is initiated by $t=0.55$. Each of the transitions between the boundary conditions is potentially problematic for a numerical scheme. When transitioning from dry to wet, if $uh$ increases much more rapidly than $h^{3/2}$ then $\tilde{B}$ in \eqref{eqn:TDP_DAE_diff} can become ill conditioned as $l^{(1)}$ and $l^{(3)}$ become essentially the same. To remedy, if either $h_R$ or $h_J$ fall below $10^{-9}$ we switch to imposing no flow, and when both increase to above $10^{-8}$ and the incident fluid motion is supercritical we switch to imposing no condition.

\begin{figure}[tp!]
	\centering
	\includegraphics{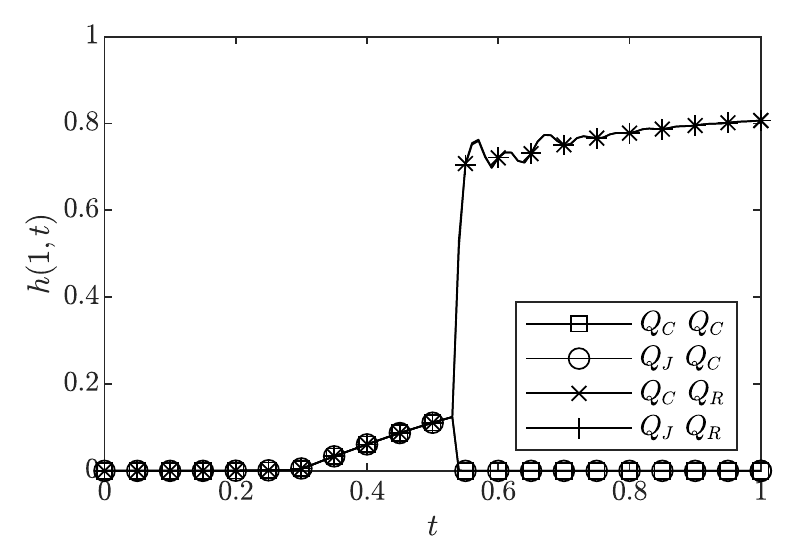}
	\caption{The depth field at the right boundary as a function of time for \cref{tp:TP_critwall}. The simulations were performed at a resolution of $J=100$ with $\mathcal{D}=1/2$, and the solution is plotted every $0.01$ time. Each solution evaluates the eigenvalues and eigenvectors at different spatial positions, in the legend the first word indicates the value used to evaluate eigenvalues for the selection of boundary conditions, the second the value used to evaluate the eigenvalues and eigenvectors for the characteristic equation.}
	\label{fig:TP_critwall_eigen}
\end{figure}

There is also the issue of evaluating eigenvalues and eigenvectors. The condition $\lambda^{(1)}|_{x=1} > 0$ on supercritical flow is a requirement that enforcing no boundary conditions is consistent with the dynamics in the bulk, so the evaluation of $\lambda^{(1)}|_{x=1}$ for this condition should be done using  $Q_{J}$. Indeed if we use $Q_R$ then the condition $uh=0$ that is used to transition from dry to wet will mean that $\lambda^{(1)}|_{x=1}<0$ and supercritical outflow is never established. For the evolution of $Q_R$ using \eqref{eqn:BC_extp_system_general} we reason that the eigenvalues and eigenvectors used in these equations should be constructed using $Q_R$. Thus eigenvalues are constructed from both $Q_R$ and $Q_{J}$. We may attempt to reduce the computational load by instead computing at the centre point $Q_C=(Q_{J}+Q_R)/2$. Simulations are plotted with different evaluation points for eigenvalues/vectors in \cref{fig:TP_critwall_eigen}, the forcing coefficient is calculated using  \cref{eqn:diffusion_outstatic} where $\whatchar{\lambda}_{R\max}$ is the largest $\abs{\whatchar{\lambda}^{(m)}}$ over all values computed for either determining the boundary condition or constructing the characteristic equations. For this problem, using $Q_R$ to construct the eigenvalues and eigenvectors used in the characteristic equation is crucial for the convergence of the scheme. For all following simulations we will do this, as well as using $Q_{J}$ for determining boundary conditions.

\begin{figure}[tp!]
	\centering
	\includegraphics{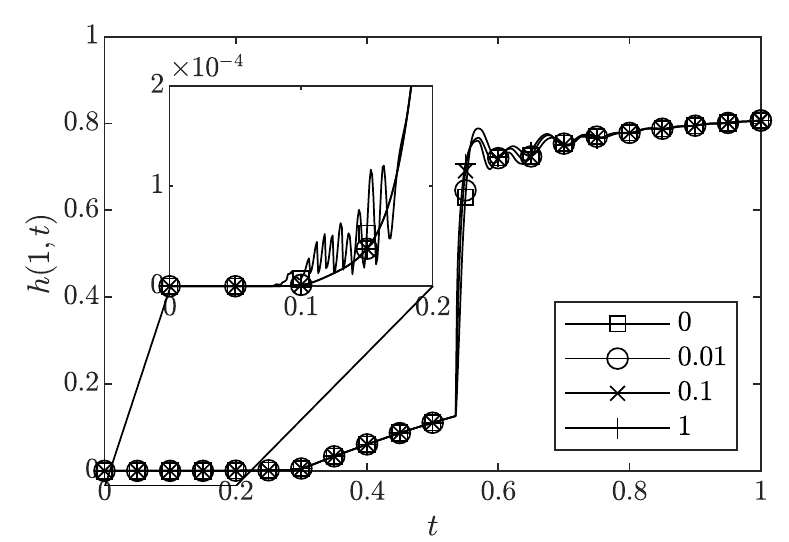}
	\includegraphics{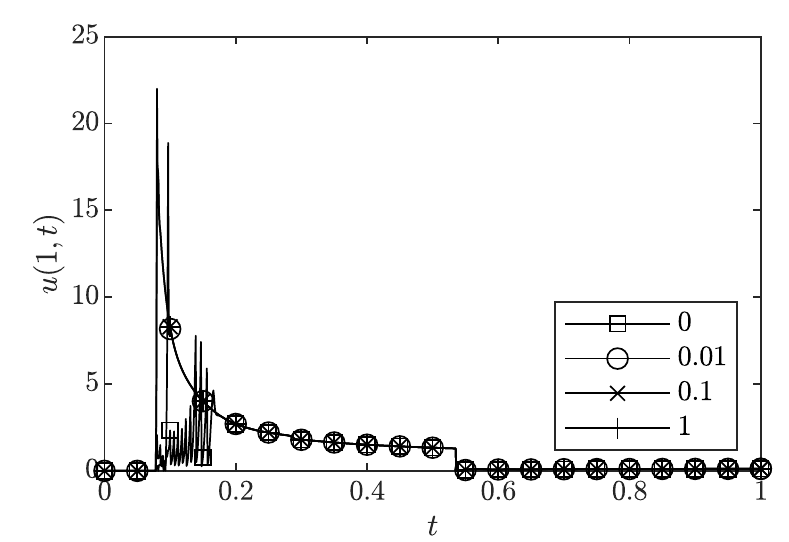}
	\caption{The depth (left) and velocity (right) fields at the right boundary as a function of time for \cref{tp:TP_critwall}. The simulations were performed at a resolution of $J=100$ with $\mathcal{D}$ as in the legend, and are plotted every $0.001$ time.}
	\label{fig:TP_critwall_diff}
\end{figure}

To investigating the influence of forcing, we plot the depth and velocity for different values of $\mathcal{D}$ in \cref{fig:TP_critwall_diff}. In the exact solution \cref{eqn:TP_critwall_exact} the fluid arrives at the barrier at time $t=1/4$, but in the simulation the fluid arrives at $t \approx 0.078$ as a result of the diffusivity of the bulk scheme. From the velocity field we see that for $\mathcal{D}=0$ there is an oscillation between no-flow and supercritical outflow conditions until $t \approx 0.18$, and for $\mathcal{D}=0.01$ this lasts until $t \approx 0.10$, while for $\mathcal{D}=0.1$ and $\mathcal{D}=1$ there is no oscillation. We conclude that forcing is necessary to stabilize transitions from subcritical to supercritical flow, keeping the boundary value relatively close to the well behaved values in the bulk. 

\begin{figure}[tp!]
	\centering
	\includegraphics{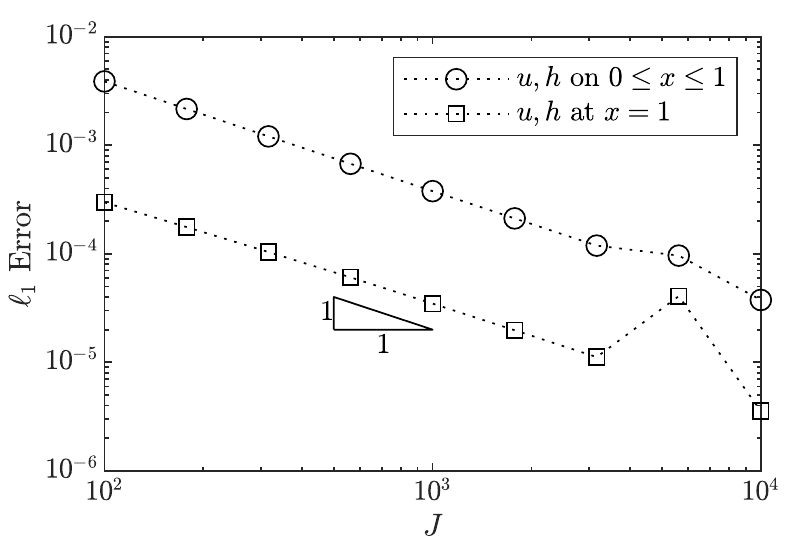}
	\caption{The $\ell_1$ error between the simulations and exact solution as a function of resolution $J$ for \cref{tp:TP_critwall} with $\mathcal{D} = 1/2$. The error is evaluated at time $t=1/2$ either across the entire domain ($0 \leq x \leq 1$) or only at the right end ($x=1$).}
	\label{fig:TP_critwall_convergence}
\end{figure}

For convergence analysis we use that \cref{eqn:TP_critwall_exact} is valid on  $0 \leq t \leq 1/2$, and \cref{fig:TP_critwall_convergence} reveals first order convergence.

\end{testproblem}

\begin{testproblem}[wall collision with particle load]	\label{tp:TP_particle}

\begin{figure}[tp!]
	\centering
	\includegraphics{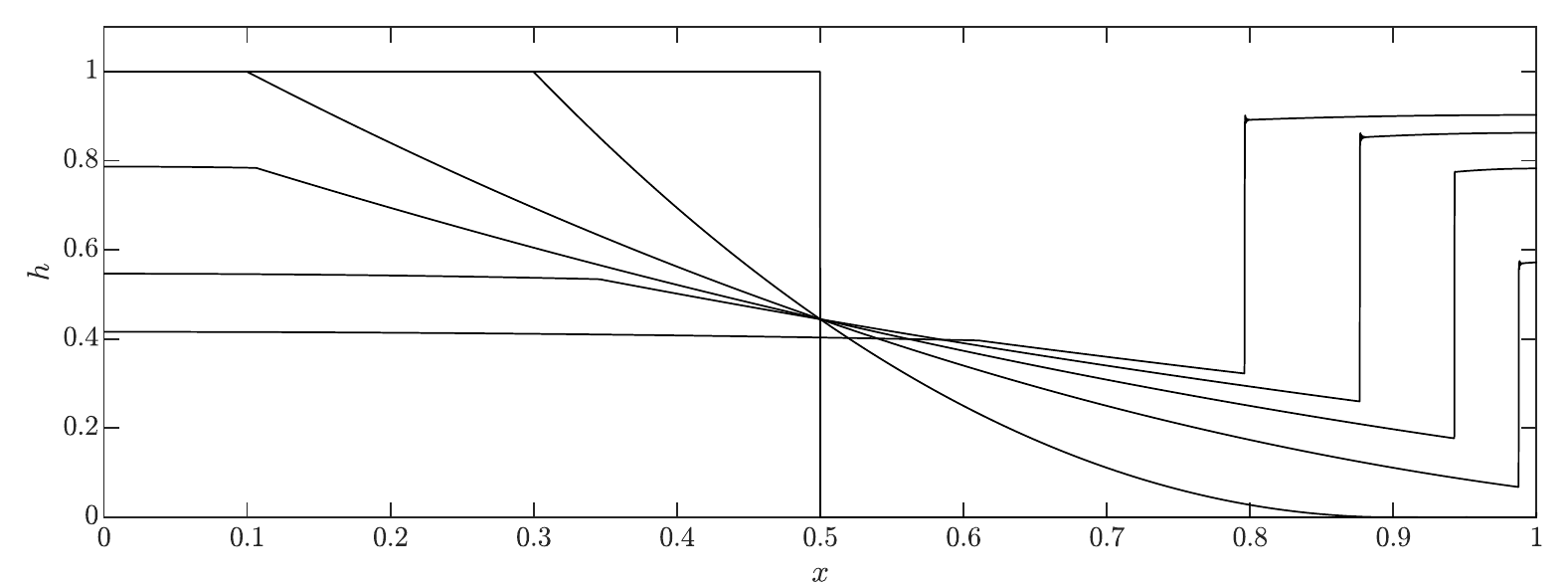}
	\caption{The depth at times $t \in \pbrc*{0, 0.2, \ldots , 1}$ for \cref{tp:TP_particle} simulated at resolution $J=10^4$ with $\mathcal{D}=1/2$. The dynamics are the same as in \cref{fig:TP_critwall} until the fluid reaches $x=1$ at time $t=1/4$, at which time the differing boundary conditions effect the dynamics.}
	\label{fig:TP_particle}
\end{figure}
\begin{figure}[tp!]
	\begin{minipage}[t]{0.5\textwidth}
		\centering
		\includegraphics{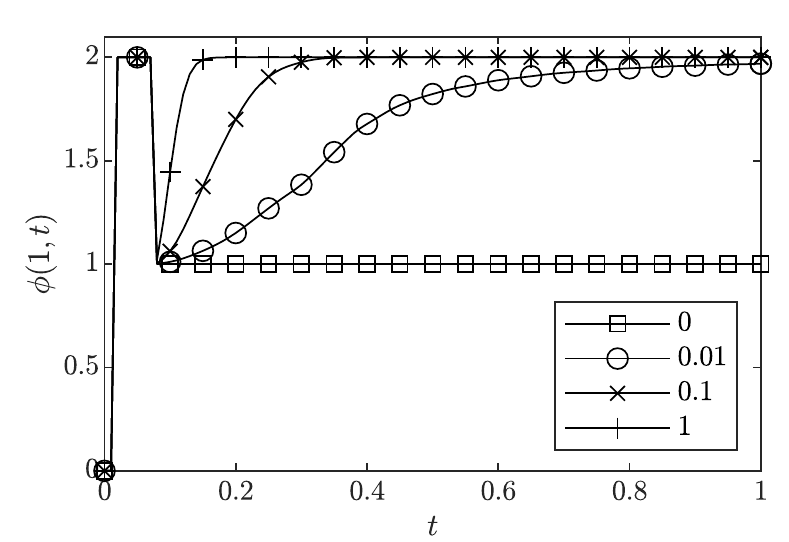}
		\captionsetup{width = 0.95\textwidth}
		\caption{The concentration field as a function of time at the right boundary for \cref{tp:TP_particle}. The simulations were performed at a resolution of $J=100$ with $\mathcal{D}$ as stated in the legend, and the solution is plotted every $0.01$ time.}
		\label{fig:TP_particle_diff}
	\end{minipage}
	\begin{minipage}[t]{0.5\textwidth}
		\centering
		\includegraphics{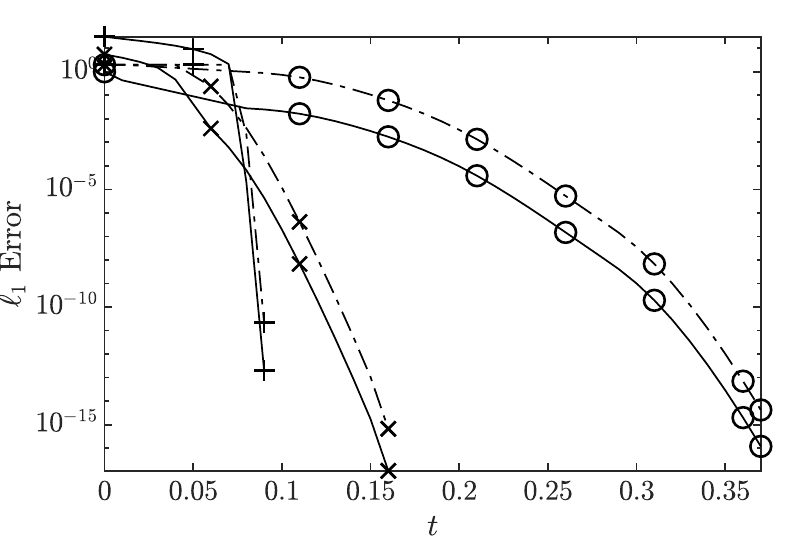}
		\captionsetup{width = 0.95\textwidth}
		\caption{The $\ell_1$ error in the concentration field $\phi$ as a function of time, plotted every $0.01$ for $0 \leq x \leq 1$, solid lines, and at $x=1$, dash dot lines. The simulations were performed with $\mathcal{D}=1/2$ and resolutions: $\bigcirc$ $J=10^2$, $\times$ $J=10^3$, $+$ $J=10.^4$. We do not plot for times when the error is precisely zero.}
		\label{fig:TP_particle_convergence}
	\end{minipage}
\end{figure}

We investigate problems encountered with the tracer field $\phi$ when it is a static characteristic. When $h>0$ then $\phi$ is constant on its characteristics, but when $h=0$ then $\phi$ is not defined. Numerically, the value of $\phi$ selected by the scheme when $h \ll 1$ is arbitrary, and we propose that forcing can be used to regularise the simulation as $h$ increases. The initial condition we consider is $h=1$, $\phi=2$ and $u=0$ for $0\leq x \leq 1/2$, and $h=0$ for $1/2<x\leq 1$, with boundary conditions $u=0$ at $x=0$ and $x=1$. The behaviour is that the fluid undergoes dam-break and collides with the barrier at time $t=1/4$, and then a backwards propagating shock is formed rapidly deepening the fluid between the shock and the barrier, as shown in \cref{fig:TP_particle} and discussed in \cite{ar_Greenspan_1978}. At the boundaries $\phi$ is constant on a static characteristic, and therefore we may expect that, for simulations without forcing, $\phi$ is constant at $x=1$ for all time. Indeed this is what is seen for $t \geq 1/4$ in \cref{fig:TP_particle_diff}, before this time the depth is very small thus $\phi$ is formally undefined and numerically unstable. Introducing even a small amount of forcing causes $\phi(1,t) \rightarrow 2$ as time passes, with $D \geq 1/10$ causing $\phi$ to be very close to $2$ by $t = 1/4$.

We plot the error as a function of time in \cref{fig:TP_particle_convergence}, and note that the error vanishes (for $\phi$ a double precision floating point) by $t=0.37$ (2 \sf), and as the resolution increases the error vanishes earlier. Because $\phi$ is not even defined analytically at $x=1$ for $t \leq 1/4$, this means that for resolutions $J \geq 10^3$ there is no error by the time $\phi$ is formally defined. Any static characteristic should converge on the correct value when forcing is included.

\end{testproblem}

\begin{testproblem}[wave in alpha field]	\label{tp:TP_alphawave}

\begin{figure}[tp!]
	\centering
	\includegraphics{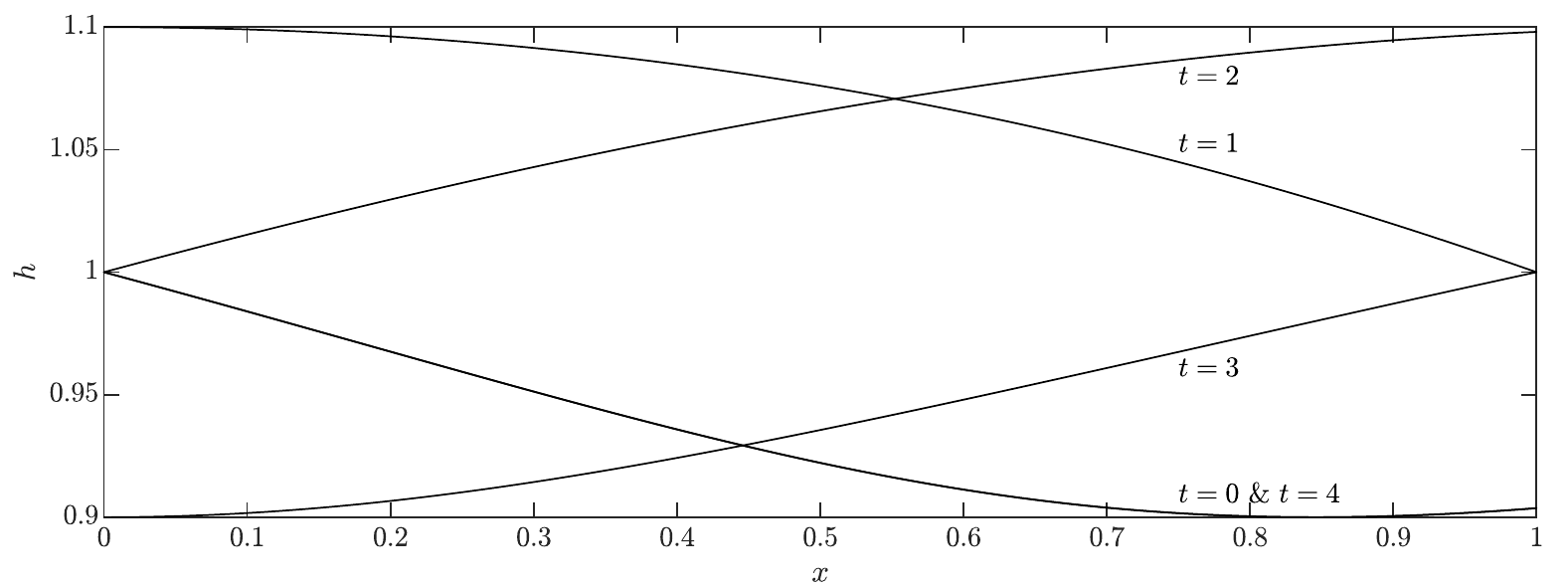}
	\caption{The depth at times $t \in \pbrc*{0,1,\ldots,4}$ (labelled) for \cref{tp:TP_alphawave} simulated at resolution $J=10^4$ with $\mathcal{D}=1/2$.}
	\label{fig:TP_alphawave}
\end{figure}
\begin{figure}[tp!]
	\centering
	\includegraphics{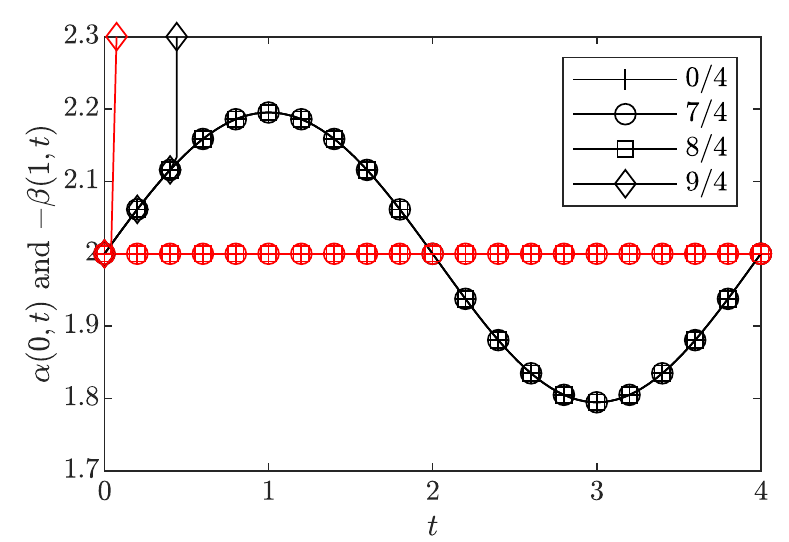}
	\includegraphics{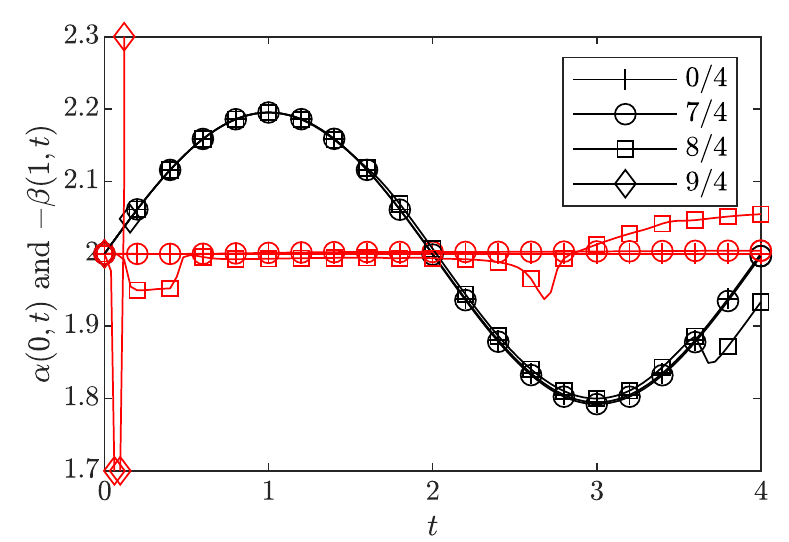}
	\caption{The inflowing invariant at either end of the domain as a function of time for \cref{tp:TP_alphawave}, in black is $\alpha(0,t)$ and in red is $-\beta(1,t)$. The simulations were performed at a resolution of $J=100$ with $\mathcal{D}$ as stated in the legend, and the solution is plotted every $0.04$ time. In the left plot \cref{eqn:diffusion_outstatic} was used, in the right \cref{eqn:diffusion_all}.}
	\label{fig:TP_alphawave_diff}
\end{figure}
\begin{figure}[tp!]
	\centering
	\includegraphics{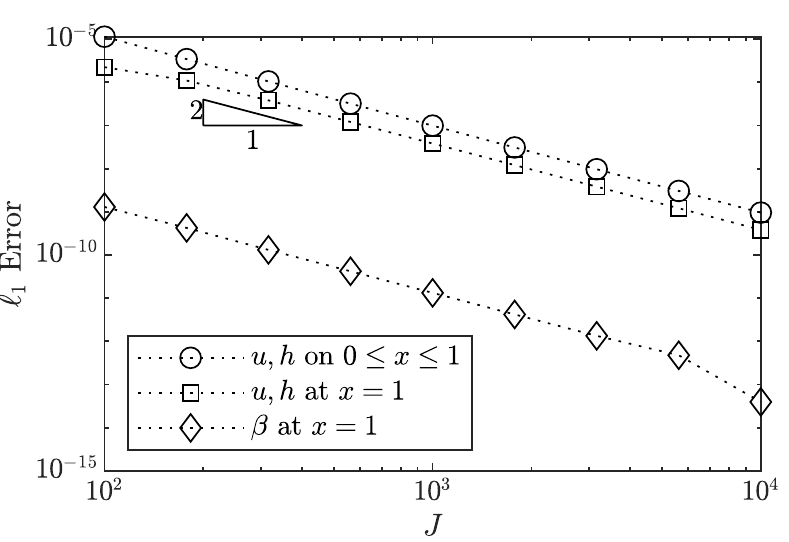}
	\includegraphics{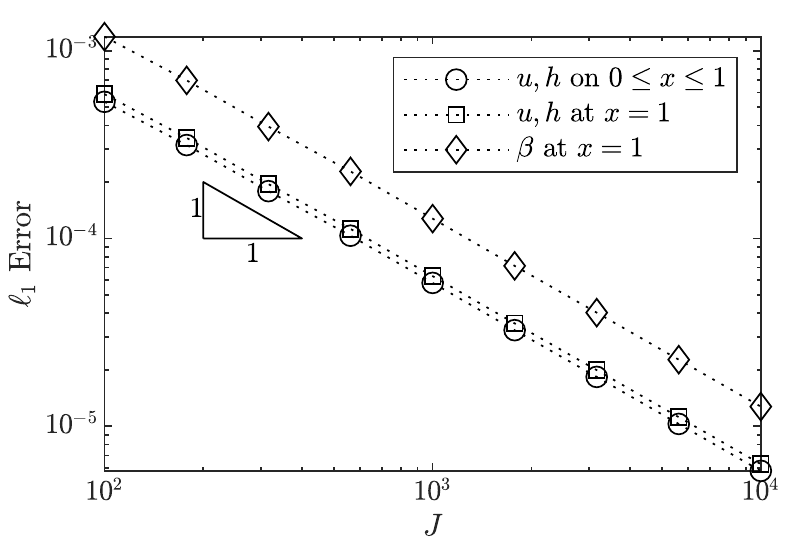}
	\caption{The $\ell_1$ error between the simulations and exact solution as a function of resolution $J$ for \cref{tp:TP_alphawave} with $\mathcal{D}=1/2$. The error is evaluated at time $t=4$ either across the entire domain ($0 \leq x \leq 1$), only at the right end ($x=1$), or only for the inflowing characteristic at the right end ($\beta(1,4)$). In the left plot \cref{eqn:diffusion_outstatic} was used, in the right \cref{eqn:diffusion_all}.}
	\label{fig:TP_alphawave_convergence}
\end{figure}

In \cref{sec:bc_discrete} we constructed bounds on the forcing coefficient for regions where the system is approximately linear. For the central-upwind scheme $C_{\max} = 1/8$ \cite{ar_Kurganov_2000,ar_Kurganov_2001}, thus we require $0 \leq \mathcal{D} < 7/4$. To verify this bound we consider a problem with smooth periodic solutions on $0 \leq x \leq 1$, imposing the boundary condition $h(0,t) = h_0(t) \eqdef 1 + \sin(\pi t/2)/10$ and requiring that $\beta=-2$ for all $x$ and $t$ (thus $u=0$ when $h=1$). Using that $u$ and $h$ are constant along $\alpha$ characteristics when $\beta$ is constant, we construct the solution \cite{ar_Lax_1957}
\begin{align}	\label{eq:TP_alphawave_soln}
	h(x,t) &= h_0(\tau),
	&
	u(x,t) &= 2\sqrt{h_0(\tau)} - 2,
	&\text{where $\tau(x,t)$ satisfies}&&
	x = (t - \tau) \pbrk*{3\sqrt{h_0(\tau)} - 2}
\end{align}
and is the time at which the $\alpha$ characteristic passing through $(x,t)$ passed through $x=0$. The simulations take the solution at $t=0$ as the initial condition, impose a non-reflecting condition at $x=1$, and evolve over a single period to $t=4$ and shown in \cref{fig:TP_alphawave}. We plot the values of the inflowing invariant at either end in \cref{fig:TP_alphawave_diff} which reveals that a $\mathcal{D}$ of $7/4$ and $8/4$ is actually stable, while for $9/4$ there is a rapidly growing instability. This gives us confidence that our bounds are sufficient for stability. We also plot for the case where forcing is included in the non-reflecting conditions using \cref{eqn:diffusion_all}, which causes a slow drift off at $7/4$, transient issues at $8/4$, and again $9/4$ is unstable. 

Examining convergence (\cref{fig:TP_alphawave_convergence}), when only including forcing in the static and outgoing characteristic fields we obtain second order convergence. We also observe that the value computed from the non-reflecting conditions, $\beta(1,t)$, has a much smaller error than the bulk does, which goes against the assumptions of \cref{sec:bc_discrete} and thus forcing is inappropriate in non-reflecting conditions. Indeed, when forcing is applied to the non-reflecting conditions we observe first order convergence, and the error in $\beta(1,t)$ is larger than that in the bulk. This is because around extrema the bulk scheme has $\order{\Delta \mathscr{y}}$ error, and when these extrema are advected through the boundary the forcing transfers this error to the boundary generating an $\order{\Delta \mathscr{y}}$ error there also, which is then advected into the domain. 

\end{testproblem}

\begin{testproblem}[lock-release] \label{tp:TP_lockrel}

\begin{figure}[tp!]
	\centering
	\includegraphics{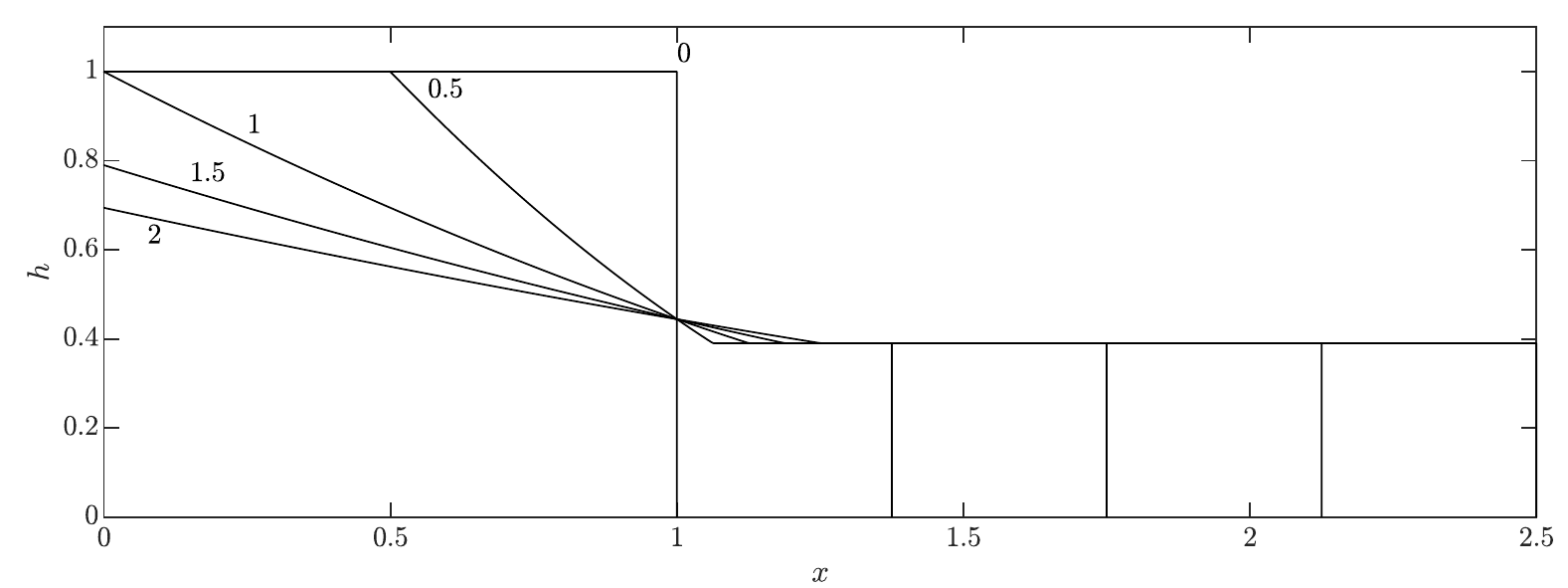}
	\caption{The depth at times $t \in \{ 0, 0.5, \ldots , 2 \}$ (labelled) for the semi-infinite \cref{tp:TP_lockrel}, simulated at $J=10^4$ with $\mathcal{D}=1/2$.}
	\label{fig:TP_lockrelinf}
\end{figure}
\begin{figure}[tp!]
	\centering
	\includegraphics{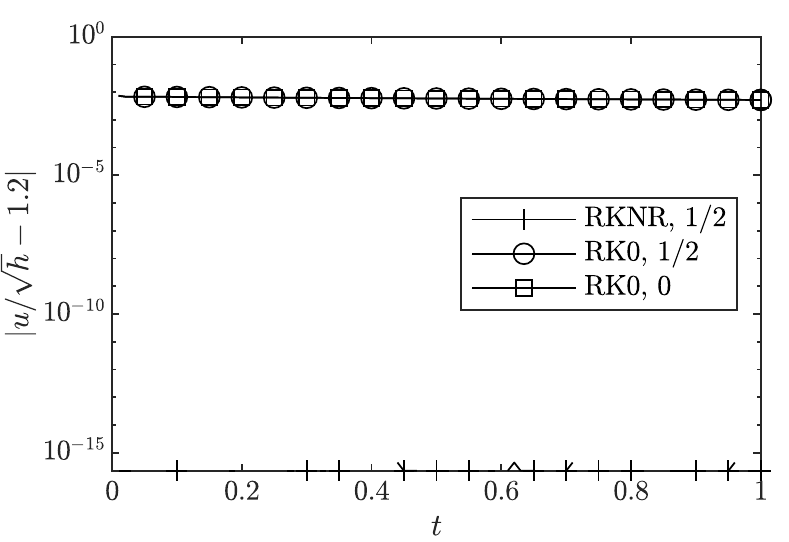}
	\caption{The error in the Froude number at $x=x_R$ as a function of time for  the semi-infinite \cref{tp:TP_lockrel}. The simulations were performed at a resolution of $J=100$, and the solution is plotted every $0.01$ time. The time-stepping algorithms and values of $\mathcal{D}$ are stated in the legend}
	\label{fig:TP_lockrel_froude}
\end{figure}
\begin{figure}[tp!]
	\centering
	\includegraphics{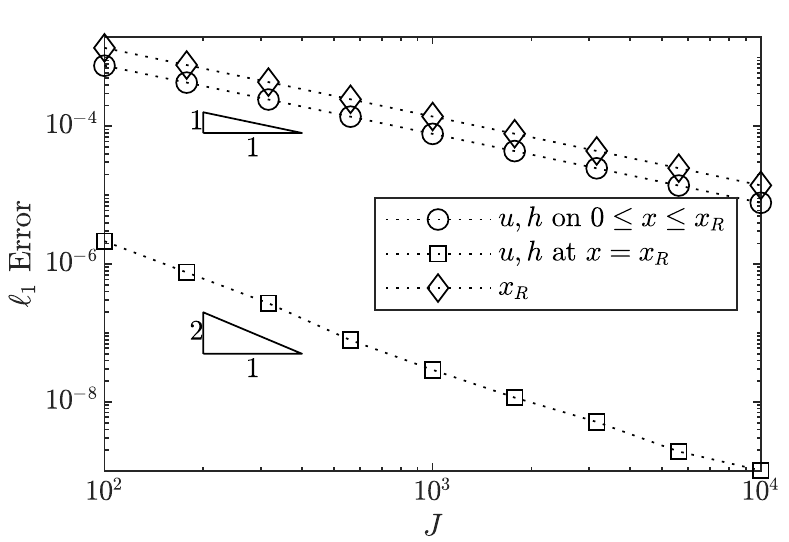}
	\includegraphics{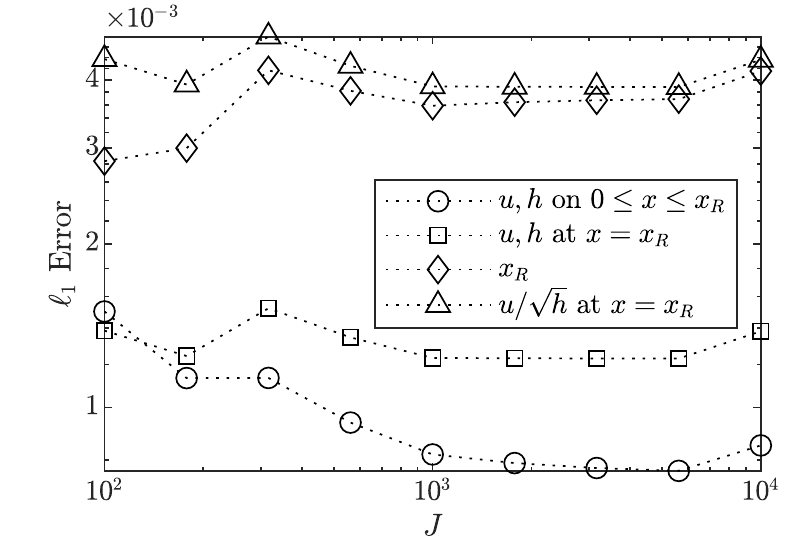}
	\caption{The $\ell_1$ error between the simulations and exact solution as a function of resolution $J$ with $\mathcal{D}=1/2$ for the semi-infinite \cref{tp:TP_lockrel}. In the left figure RKNR time-stepping was used, in the right RK0. The error is evaluated at time $t=2$ for the expressions in the legends.}
	\label{fig:TP_lockrelinf_convergence}
\end{figure}

Lastly, we simulate the lock-release of a gravity current. We begin with the semi-infinite configuration (\cref{fig:TP_lockrelinf}), which we simulate on an evolving domain $0 \leq x \leq x_R(t)$, with initial conditions $h(x,0)=1$, $u(x,0)=0$, $x_R(0) = 1$, non-reflecting conditions at $x=0$, and $u = \dot{x}_R = Fr \sqrt{h}$ at $x=x_R$ (see \cite{ar_Benjamin_1968}); we take $Fr = 1.2$. From \cref{fig:TP_lockrel_froude} we see that, without projection, the boundary condition drifts over a short time interval after initialisation. While the error generated is small, the local invariance of the system around $(x,t)=(1,0)$ with respect to equal rescaling of $x$ and $t$ indicates that this error will not reduce with resolution. This is what is observed in \cref{fig:TP_lockrelinf_convergence}, with projection the solution converges at first order, while without the solution \emph{does not converge}.

\begin{figure}[tp!]
	\centering
	\includegraphics{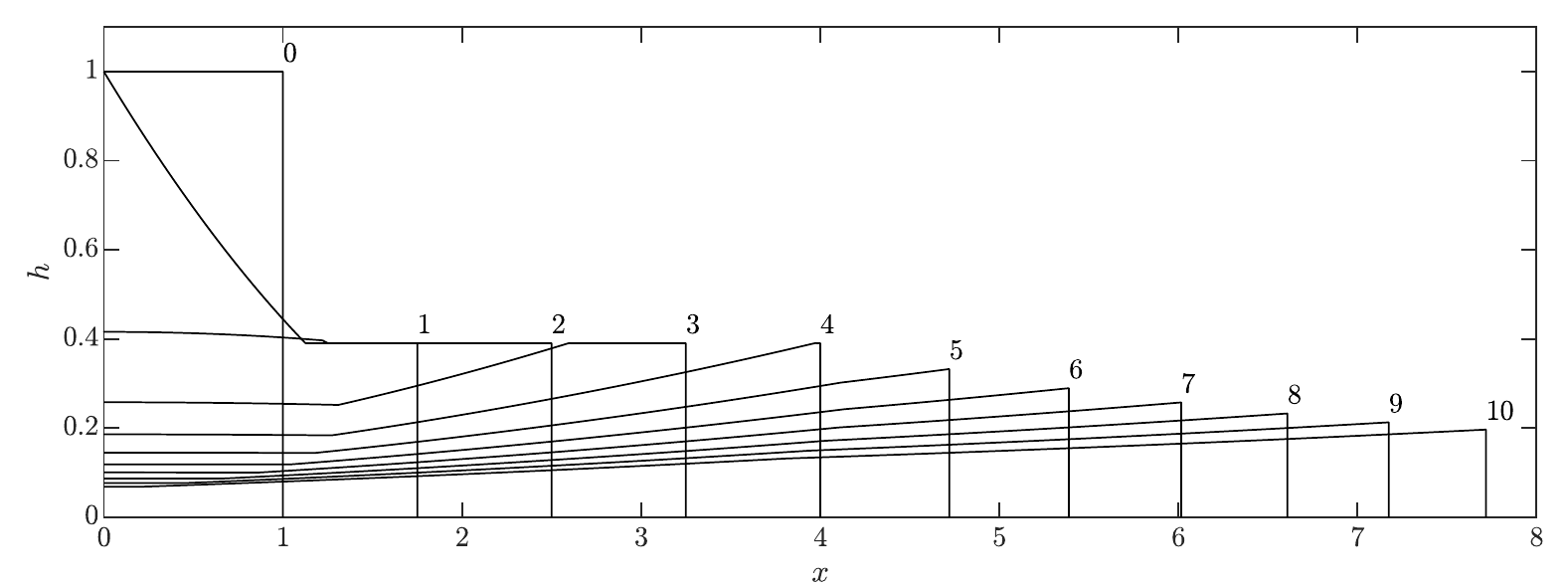}
	\caption{The depth at times $t \in \{ 0, 1, \ldots , 10 \}$  (labelled) for the finite \cref{tp:TP_lockrel}, simulated at $J=10^4$ with $\mathcal{D}=1/2$.}
	\label{fig:TP_lockrel}
\end{figure}
\begin{figure}[tp!]
	\centering
	\includegraphics{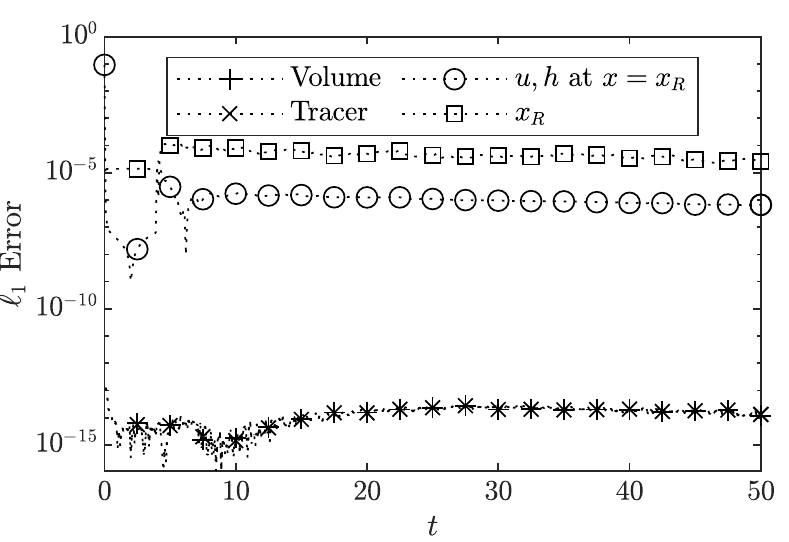}
	\caption{The error as a function of time for the finite \cref{tp:TP_lockrel}, simulated at resolution $J=10^4$ with $\mathcal{D}=1/2$. The error is evaluated at time intervals of $0.1$ for the total volume of fluid in the domain (Volume), the total amount of tracer in the domain (Tracer), $u$ and $h$ at the right end ($x=x_R$), and the location of the front ($x_R$). Comparison data for the evolution at $x=x_R$ courtesy of A. J. Hogg, see \cite{ar_Hogg_2006}.}
	\label{fig:TP_lockrel_error}
\end{figure}

To demonstrate the effectiveness of our scheme, we present a simulation of the lock-release problem on a finite domain (\cref{fig:TP_lockrel}); that is we have initial conditions $h(x,0)=\phi(x,0)=1$, $u(x,0)=0$, $x_R(0) = 1$, and boundary conditions $u=0$ at $x=0$ and $u = \dot{x}_R = Fr \sqrt{h}$ at $x=x_R$. This problem was solved exactly in \cite{ar_Hogg_2006}, revealing that at time $t = 4.05$ (3 \sf) a gradient discontinuity collides with $x=x_R$ causing the boundary value to begins evolving. A small amount of error is generated at the front during this transition, see \cref{fig:TP_lockrel_error}, but it is three orders of magnitude lower than the error seen in \cref{fig:TP_lockrelinf_convergence} for RK0 at all resolutions, and the boundary condition is resolved to machine precision. The error in $x_R$ caused by transient issues around the initial dam-break, and the transition from constant to evolving height at the front. In \ref{fig:TP_lockrel_error} we also plot the error in the total volume and tracer load in the system, and show that these quantities are conserved to machine precision, which is due to the transformation performed in \cref{sec:SW}.

\end{testproblem}

%% file: Sections/Discussion.tex
\section{Summary}	\label{sec:Discussion}

We have presented two improvements to the method of \cite{ar_Thompson_1987} for the implementation of non-linear algebraic and differential-algebraic boundary conditions for hyperbolic systems.

Extrapolation and forcing was included in the characteristic equations used to complete the system at the boundary, and in \cref{sec:bc_discrete} we carefully constructed discretized expressions for the derivatives that ensure the boundary value is consistent with the solution in the bulk. In \cref{tp:TP_critwall,tp:TP_particle}, forcing in the boundary conditions was necessary to properly resolve the dynamics in simulations of the shallow water equations. For \cref{tp:TP_critwall} we showed that without forcing the simulations could not transition from subcritical to supercritical flow, but with it they could. For \cref{tp:TP_particle} the value of the static characteristic was found to rapidly converge with forcing, while not including forcing meant that the the solution did not converge. We presented indicative stability analysis yielding bounds on the coefficient \eqref{eqn:BCSP_onesided_convergence_weakened}, and we confirmed that these bounds did indeed yield stable simulations in \cref{tp:TP_alphawave}. It was also found that, while including forcing in the static and outgoing characteristic equations greatly improves the properties of the simulation, including it in the non-reflecting conditions of incoming characteristics lead to a loss of accuracy, thus we advise setting the coefficient to zero in this case. No special properties of the shallow water equations were required to achieve convergence, and for this reason we expect that the method presented is able to produce convergent simulations for a wide class of hyperbolic systems.

Projection methods were discussed in \cref{sec:TDP}, including the exposition of the RKNR scheme originally presented in \cite{ar_Skevington_F001_Draining} and included here in greater detail. From our investigation we concluded that, for second order schemes, our RKNR scheme was best (implicitly solving the algebraic system during each sub-step), and for third order schemes RKP1 was best (projecting onto the manifold after completion of each step). In test problem \ref{tp:TP_lockrel} we established that projection is necessary to produce convergent simulations because using a standard Runge-Kutta scheme yielded $\order{1}$ errors. While the error may be small, for problems that are dependant on a delicate balance in the boundary conditions, such as in \cite{ar_Skevington_F001_Draining}, the presence of any error in the evaluation of the boundary condition can dominate the computation of the solution and so a projection scheme must be used.

%% file: Sections/Acknowledgements.tex
\section*{Acknowledgements}

This work was supported by the EPSRC (grant number EP/M506473/1). The author would also like to thank A. J. Hogg and C. G. Johnson for their constructive comments regarding drafts of this article.

%% file: Sections/RK_contraction.tex
\section{Proof of \cref{thm:FV_time_step_stability_full}} \label{sec:RK_contraction_proof}

We begin by establishing a general result for the behaviour of bounded semi-norms under the RK scheme \cref{eqn:FV_RK}.

\begin{lemma}\label{thm:FV_time_step_stability_part}
\begin{subequations}
	Let $z^{n[s]}$ and $\hat{z}(t)$ be as in \cref{thm:FV_time_step_stability_full}, and
	\begin{align}	\label{eqn:FV_thm_timestepstabilitypart_boundfunc}
		\tilde{f}(K;s;a,b,c)
		&= a K^s \prod_{\beta = 0}^{s-1} \ppar*{1-\eta^{[\beta]}} + b \sum_{\alpha = 0}^{s-1} K^{\alpha} \eta^{[s-\alpha-1]} \pbrk*{\prod_{\beta=s-\alpha}^{s-1} \ppar*{1-\eta^{[\beta]}}} + c \sum_{\alpha = 0}^{s-1} K^{\alpha} \pbrk*{\prod_{\beta=s-\alpha}^{s-1} \ppar*{1-\eta^{[\beta]}}}
	\end{align}
	where $\tilde{f}(K;0;a,b,c) = a$. If, for some constants $K,\tilde{K} \geq 0$,
	\begin{equation}\label{eqn:FV_thm_timestepstabilitypart_step}
		\norm*{\mathscr{E}(z^{n[s]}) - \hat{z}\ppar*{t^{n[s]} + \Delta t^n}} \leq K \norm*{z^{n[s]} - \hat{z}\ppar*{t^{n[s]}}} + \tilde{K} 
	\end{equation}
	for all $s \in \{0,1,\ldots,S-1\}$, then for all $s \in \{0,1,\ldots,S\}$
	\begin{equation}\label{eqn:FV_thm_timestepstabilitypart_result}
		\norm*{z^{n[s]} - \hat{z}\ppar*{t^{n[s]}}} \leq \tilde{f}(K;s;1,1,0)\norm*{z^{n[0]} - \hat{z}\ppar*{t^{n[0]}}} + \tilde{f}(K;s;0,-1,1) \tilde{K}  + \tilde{f}\ppar*{K;s;0,-\textfrac{1}{2},1} L_z \ppar*{\Delta t^n}^2.
	\end{equation}
\end{subequations}\end{lemma}
\begin{proof}
	We start by bounding after a sub-step.
	\begin{gather*}
		\begin{split}
			z^{n[s+1]} - \hat{z}\ppar*{t^{n[s+1]}}
			&= 
			\ppar*{1-\eta^{[s]}} \pbrk*{ \mathscr{E}\ppar*{z^{n[s]}} - \hat{z}\ppar*{t^{n[s]} + \Delta t^n} }
			+ \eta^{[s]} \pbrk*{z^{n [0]} - \hat{z}\ppar*{t^{n [0]}} }
			\\&\qquad\qquad
			+ \ppar*{1-\eta^{[s]}}\hat{z}\ppar*{t^{n[s]} + \Delta t^n} - \hat{z}\ppar*{t^{n[s+1]}} +\eta^{[s]} \hat{z}\ppar*{t^{n[0]}}.
		\end{split}
	\shortintertext{Using that}
		\norm*{\ppar*{1-\eta^{[s]}}\hat{z}\ppar*{t^{n[s]} + \Delta t^n} - \hat{z}\ppar*{t^{n[s+1]}} +\eta^{[s]} \hat{z}\ppar*{t^{n[0]}}}
		 \leq \ppar*{ 1 - \textfrac{1}{2} \eta^{[s]}} L_z \ppar*{\Delta t^n}^2
	\end{gather*}
	where $\hat{z}$ was Taylor expanded around $t^{n[0]}$, we establish
	\begin{align*}
		\norm*{z^{n[s+1]} - \hat{z}\ppar*{t^{n[s+1]}}} \leq \ppar*{1-\eta^{[s]}} K \norm*{ z^{n[s]} - \hat{z}\ppar*{t^{n[s]}} } + \eta^{[s]} \norm*{ z^{n [0]} - \hat{z}\ppar*{t^{n [0]}} } + \ppar*{1-\eta^{[s]}}\tilde{K} + \ppar*{1-\textfrac{1}{2}\eta^{[s]}}L_z \ppar*{\Delta t^n}^2
	\end{align*}
	This gives a recurrence relation, the bound after a sub-step expressed in terms of the values prior to the sub-step. We next use induction to show \cref{eqn:FV_thm_timestepstabilitypart_result}. Examining the base case $s=0$, we have that $\tilde{f}(K;0;a,b,c)=a$, thus our proposed expression is trivially true. Performing the inductive step from $s$ to $s+1$
	\begin{align*}\begin{split}
		\norm*{z^{n[s+1]} - \hat{z}\ppar*{t^{n[s+1]}}} \leq 
		& \pbrk*{\ppar*{1-\eta^{[s]}} K \tilde{f}(K;s;1,1,0) + \eta^{[s]}}\norm*{z^{n[0]} - \hat{z}\ppar*{t^{n[0]}}} \\
		& + \pbrk*{\ppar*{1-\eta^{[s]}} K \tilde{f}(K;s;0,-1,1) + 1 - \eta^{[s]}} \tilde{K}  \\
		& + \pbrk*{\ppar*{1-\eta^{[s]}} K \tilde{f}\ppar*{K;s;0,-\textfrac{1}{2},1} + 1 - \textfrac{1}{2} \eta^{[s]}} L_z \ppar*{\Delta t^n}^2.
	\end{split}\end{align*}
	Examining the coefficients
	\begin{align*}\begin{split}
		&\ppar*{1-\eta^{[s]}} K \tilde{f}(K;s;a,b,c) + b \eta^{[s]} + c	\\
		&= a K^{s+1} \prod_{\beta = 0}^{s} \ppar*{1-\eta^{[\beta]}}
		+ b \ppar*{ \ppar*{1-\eta^{[s]}} K \sum_{\alpha = 0}^{s-1} K^{\alpha} \eta^{[s-\alpha-1]} \pbrk*{\prod_{\beta=s-\alpha}^{s-1} \ppar*{1-\eta^{[\beta]}}} + \eta^{[s]}}
		+ c \ppar*{ \ppar*{1-\eta^{[s]}} K \sum_{\alpha = 0}^{s-1} K^{\alpha} \pbrk*{\prod_{\beta=s-\alpha}^{s-1} \ppar*{1-\eta^{[\beta]}}} + 1}	\\
		&= a K^{s+1} \prod_{\beta = 0}^{s} \ppar*{1-\eta^{[\beta]}}
		+ b \sum_{\alpha = 0}^{s} K^{\alpha} \eta^{[s-\alpha]} \pbrk*{\prod_{\beta=s-\alpha+1}^{s} \ppar*{1-\eta^{[\beta]}}}
		+ c \sum_{\alpha = 0}^{s} K^{\alpha} \pbrk*{\prod_{\beta=s+1-\alpha}^{s} \ppar*{1-\eta^{[\beta]}}}
		= \tilde{f}(K;s+1;a,b,c).
	\end{split}\end{align*}
	This concludes the proof.
\end{proof}

To prove \cref{thm:FV_time_step_stability_full} we examine \cref{thm:FV_time_step_stability_part} for the case $s=S$ for the schemes in \cref{tab:RKcoeff}, for which
\begin{subequations}\begin{align}
	\tilde{f}(K;S;a,b,c) &= aK + c
	&\text{for}&&
	S &= 1,	\\
	\tilde{f}(K;S;a,b,c) &= \textfrac{1}{2}aK^2 + \textfrac{1}{2} c K + \textfrac{1}{2}b + c
	&\text{for}&&
	S &= 2,	\\
	\tilde{f}(K;S;a,b,c) &= \textfrac{1}{6}aK^3 + \textfrac{1}{6} c K^2 + \ppar*{\textfrac{1}{2}b + \textfrac{2}{3}c}K + \textfrac{1}{3}b + c
	&\text{for}&&
	S &= 3.
\end{align}\end{subequations}
For all $\tilde{f}(1;S;1,1,0) = 1$ and $\tilde{f}(1;S;0,-1,1) = 1$, while $\tilde{f}(1;1;0,-1/2,1) = 1$, $\tilde{f}(1;2;0,-1/2,1) = 5/4$, and $\tilde{f}(1;3;0,-1/2,1) = 17/12$.

%% file: Sections/Implementation.tex
\section{Implementation}	\label{sec:IMP}

Here we discuss our implementation of the scheme using second order RKNR time-stepping (\cref{sec:TDP_RKNR}) and the central-upwind scheme \cite{ar_Kurganov_2001} on a uniform grid. We use the word \emph{mode} to refer to the choice of a particular piece of a piecewise defined algebraic boundary condition, and we refer to $\whatchar{h}$ and $\widehat{\phi h}$ as \emph{positivity preserved} variables, if they should be set to zero we call them a \emph{zeroed} variable. The algebraic conditions are transformed by (\cref{sec:TED})
\begin{equation}
	\mathscr{g}(\mathscr{Q}_R , \mathscr{Q}_{\mathscr{y}R} , \dot{x}_R , x_R , t) = g(\mathscr{T}^{-1} \mathscr{Q}_R , \mathscr{T}^{-1} \mathscr{Q}_{\mathscr{y}R}/\mathscr{l} , \dot{x}_R , x_R , t).
\end{equation}

\paragraph{Initialisation}

To initialise the values in the bulk we average the initial conditions over each cell. To initialise the domain endpoint and then solve the optimisation problem (initialised with $\mathscr{Q}_{R0}$ from initial conditions and $\dot{x}_R=0$)
\begin{align}	\label{eqn:IMP_minimisation}
	&\text{minimise}&&
	(h_R-h_{R0})^2 + (\phi_R-\phi_{R0})^2 + (u_R-u_{R0})^2&
	&\text{subject to}&
	\mathscr{g}&=0,&
	\whatchar{h}_R& \geq 10^{-8},&
	\widehat{\phi h}_R& \geq 10^{-8}.
\end{align}
If the boundary condition is defined piecewise with conditions, we iterate solving the optimisation problem upon each piece until the conditions are satisfied. After minimising, if $\whatchar{h}_R = 10^{-8}$ then we zero it, similarly for $\widehat{\phi h}_R$.

\paragraph{Single time-step}

To perform an Euler time-step we begin by time-stepping the boundary conditions. To do this we first construct the full set of characteristic equations and non-reflecting conditions, \ie we construct $\tilde{B}$ and $\tilde{b}$ from \cref{eqn:TDP_DAE_diff} as though we were constructing non-reflecting conditions \cref{eqn:BC_extp_system_reflect}. For the purposes of this construction only, we make the transformation $\whatchar{h} \mapsto \max(\whatchar{h},\delta)$ and $\widehat{\phi h} \mapsto \max(\widehat{\phi h},\delta)$ where $\delta = 10^{-8}$, which yields distinct eigenvalues and linearly independent eigenvectors. We then iterate the following until the choice of boundary condition and zeroed values has stabilised, and the change in $\Delta v$ is below $10^{-12}$:
\begin{enumerate*}[label=(\arabic*)]
	\item compute stepped values
	\item check to see if any positivity preserved variables are less than $10^{-8}$, if so then they become zeroed
	\item change mode if required
	\item assemble system to solve
	\item if any values are zeroed then overwrite fastest outgoing characteristic equations with condition that these values should be zero
	\item compute next value of $\Delta v$
\end{enumerate*}
. Any zeroed values are then set to zero, and if $\whatchar{h}_R > 10^4 \whatchar{h}_{J}$ then we set $\whatchar{h}_R = 10^4 \whatchar{h}_{J}$ (and the same for $\widehat{\phi h}_R$) to prevent the endpoint value being many orders of magnitude greater than the closest bulk value, which has been found to cause problems. 

We then compute the reconstruction using the generalised minmod reconstruction recommended in \cite{ar_Kurganov_2001}, the gradient in each cell being
\begin{align}
	[\mathscr{Q}_\mathscr{y}]_j  = \minmod\ppar*{ \theta \frac{\mathscr{Q}_j - \mathscr{Q}_{j-1}}{\Delta \mathscr{y}} , \frac{\mathscr{Q}_{j+1} - \mathscr{Q}_{j-1}}{2 \Delta \mathscr{y}} , \theta \frac{\mathscr{Q}_{j+1} - \mathscr{Q}_j}{\Delta \mathscr{y}} },
\end{align}
where we take $\theta = 3/2$. This requires values in the cells on either side, but the cells adjacent to the boundaries do not have both neighbours. To fix this we extrapolate a value beyond the domain end (this extrapolation is different to the one used in \cref{sec:spatial}) using linear interpolation from the values at $\mathscr{y} \in \{ \mathscr{y}_J , 1\}$ to produce a value at the point $\Delta \mathscr{y}/2$ beyond the domain end. The fluxes are then evaluated using the central-upwind scheme, except for at the domain ends at which we have the value of $\mathscr{Q}$ and so evaluate the flux directly. The sources are evaluated for each cell using the average value (i.e. the integral average of the source across each cell is evaluated using the midpoint method), and from this the time derivative for the bulk points is deduced and the time-step made. Finally, for any cell with $\whatchar{h} < 10^{-8}$ we decrease $\widehat{uh}$ by $10 \%$ after every time-step to prevent shallow cells gaining unbounded velocities.

\paragraph{Time Evolution}

To evolve through time we simply perform many Runge-Kutta time-steps sequentially. To ensure that the CFL condition is always met, we attempt a time-step at $95\%$ the maximum size reported by any Euler step in the previous RK step. If any Euler step reports a smaller value of the maximum time-step size than the one used then we reduce it accordingly and reattempt. When we approach a time at which we wish to store an output then we step onto this time, so long as such a step would not be larger than $95\%$ the maximum step size.

When boundary conditions are used that can change the number of conditions applied, or contain a discontinuity in the manifold, then we use event functions to detect the time at which the change occurs, an event function changing sign at the change in boundary condition. If the normal Runge-Kutta time-step would result in a change of sign of the event function then instead we perform an Euler time-step, use root finding to find the time at which the event is triggered, and take a step of this size. We then solve the minimisation problem \eqref{eqn:IMP_minimisation} to project onto the new boundary conditions, initialising with the values reached in the simulation, and then continue evolution.

\paragraph{Energy boundary condition}

For the energetic outflow condition in \cref{tp:TP_critwall} the boundary condition $\dot{x}_R = 0$ is always enforced. Besides this we employ the energy
\begin{equation}
	\widehat{\Delta E} = \widehat{uh} + 2 \mathscr{l} \whatchar{h}^2 \ppar*{ \whatchar{h} - \mathscr{l} h_b } - 3 \mathscr{l}^{2/3} \whatchar{h}^2 \widehat{uh}^{2/3}
\end{equation}
which is equal to $2 \mathscr{l}^4 h^2 \Delta E$. Supercritical outflow must be enforced when $\lambda^{(1)} > 0$, $\Delta E > 0$ and the depth is above some small threshold. During evolution we transition into and out of the supercritical outflow by using the event detection functions

\begin{subequations}\begin{align}
	-\min\ppar*{ \frac{\whatchar{\lambda}_J^{(1)}}{\min\ppar*{\sqrt{h_R},10^{-4}}} , \widehat{\Delta E_J} \cdot 10^{12} , \min(\whatchar{h}_J,\whatchar{h}_R) \cdot 10^{6} - 1 } + 10^{-3}
\shortintertext{and}
	\min\ppar*{ \frac{\whatchar{\lambda}_J^{(1)}}{\min\ppar*{\sqrt{h_R},10^{-4}}} , \widehat{\Delta E_J} \cdot 10^{12} , \min(\whatchar{h}_J,\whatchar{h}_R) \cdot 10^{7} - 1 } + 10^{-3}
\end{align}\end{subequations}
respectively. Note that $\whatchar{\lambda}_J^{(1)}$ is replaced by $\whatchar{\lambda}_C^{(1)}$ in some simulations.

When the outflow is subcritical we either impose $u=0$ (for $h<h_b$) or $\Delta E = 0$. The solution to $\Delta E = 0$ must be in the range $0 \leq u \leq \sqrt{h}$, which we enforce by extrapolating the function for values outside of this range, \ie we impose $\mathscr{g}=0$ where
\begin{align}
	\mathscr{g} &=
	\begin{cases}
		\widehat{uh} 
		& \text{for } \whatchar{h} < \mathscr{l} h_b,
		\\
		\eval{ \widehat{\Delta E} }_{\widehat{uh} = \delta} + \eval{ \pdv{\widehat{\Delta E}}{\widehat{uh}} }_{\widehat{uh} = \delta} \cdot (\widehat{uh} - \delta) 
		& \text{for } \widehat{uh} < \delta,
		\\
		\eval{ \widehat{\Delta E} }_{\widehat{uh} = \mathscr{l}^{1/2} \whatchar{h}^{3/2} - \delta} +  \eval{ \pdv{\widehat{\Delta E}}{\widehat{uh}} }_{\widehat{uh} = \mathscr{l}^{1/2} \whatchar{h}^{3/2} - \delta} \cdot (\widehat{uh} - \mathscr{l}^{1/2} \whatchar{h}^{3/2} + \delta) 
		& \text{for } \widehat{uh} >\mathscr{l}^{1/2} \whatchar{h}^{3/2} - \delta,
		\\
		\widehat{\Delta E}
		& \text{otherwise}
	\end{cases}
\end{align}
with $\delta = 10^{-8}$.